\documentclass{article}
\usepackage{iclr2027_conference,times}

\usepackage{algorithm}
\usepackage{algorithmic}
\usepackage{graphicx}

\usepackage{hyperref}
\usepackage{url}
\usepackage{amsmath,amsthm,amsfonts}
\usepackage{makecell}
\usepackage{subcaption}
\usepackage{paralist}
\usepackage{pifont}
\usepackage[table]{xcolor}
\usepackage{threeparttable}
\usepackage{multirow,booktabs}

\theoremstyle{plain}
\newtheorem{theorem}{Theorem}

\newtheorem{lemma}{Lemma}
\newtheorem{corollary}{Corollary}
\theoremstyle{definition}
\newtheorem{definition}{Definition}
\newtheorem{assumption}{Assumption}
\theoremstyle{remark}

\usepackage{amsmath,amsfonts,bm}

\def\eqref#1{equation~\ref{#1}}

\def\1{\bm{1}}

\DeclareMathAlphabet{\mathsfit}{\encodingdefault}{\sfdefault}{m}{sl}
\SetMathAlphabet{\mathsfit}{bold}{\encodingdefault}{\sfdefault}{bx}{n}

\newcommand{\E}{\mathbb{E}}

\newcommand{\R}{\mathbb{R}}

\DeclareMathOperator*{\argmax}{arg\,max}
\DeclareMathOperator*{\argmin}{arg\,min}

\def \x {\mathbf{x}}
\def \g {\mathbf{g}}
\def \bg {\mathbf{\bar{g}}}
\def \z {\mathbf{z}}
\def \u {\mathbf{u}}
\def \w {\mathbf{w}}
\def \P {\mathcal{P}}
\def \F {\mathcal{F}}
\def \G {\mathcal{G}}
\def \W {\mathcal{W}}
\def \A {\mathcal{A}}
\def \q {\mathbf{q}}
\def \s {\mathbf{s}}
\def \C {\mathcal{C}}
\def \V {\mathcal{V}}
\def \wb {\bar{\w}}
\def \qb {\bar{\q}}
\def \e {\mathbf{e}}
\def \Lb {\mathbf{L}}
\def \hLb {\mathbf{\hat{L}}}
\def \bLb {\mathbf{\bar{L}}}
\def \epswdist  {\varepsilon_{w,t}}
\def \epsqdist  {\varepsilon_{q,t}}
\def \Fb {\mathbf{F}}
\def \hFb {\mathbf{\hat{F}}}
\def \bFb {\mathbf{\bar{F}}}
\def \Xb {\mathbf{X}}
\def \hXb {\mathbf{\hat{X}}}
\def \bXb {\mathbf{\bar{X}}}
\def \Regq {\mathrm{Reg}_{q,T}}
\def \Regw {\mathrm{Reg}_{w,T}}
\newcommand\inner[2]{\langle #1, #2 \rangle}
\def \hats {\hat{s}}
\def \tildesb {\tilde{\s}}

\title{Group Distributionally Robust Optimization over Decentralized Networks}

\author{
Haomin Bai$^{1,2}$ \quad
Junyan Sun$^{1,2}$ \quad
Quan Cheng$^{1,2}$ \quad
Yuanyu Wan$^{3}$ \quad
Lijun Zhang$^{1,2}$\thanks{Lijun Zhang is the corresponding author.}\\
$^1$State Key Laboratory of Novel Software Technology, Nanjing University \\
$^2$School of Artificial Intelligence, Nanjing University \\
$^3$School of Software Technology, Zhejiang University \\
\texttt{\{baihm,sunjy,chengq,zhanglj\}@lamda.nju.edu.cn}, \texttt{wanyy@zju.edu.cn}
}

\iclrfinalcopy
\begin{document}

\maketitle

\begin{abstract}
We study group distributionally robust optimization (GDRO) over decentralized networks, where workers collaborate to learn robust models by minimizing the maximum local risk.
Existing methods optimize only dual-regularized surrogates and establish guarantees for a globally averaged model that is not available to workers.
In this paper, we directly optimize the original GDRO objective and establish guarantees for local models in semi-decentralized (SemiDec) and fully decentralized (Dec) settings.
A key difficulty is that applying the primal-dual decomposition from centralized GDRO yields worker-dependent losses whose evaluation requires unavailable cross-worker information.
To address this challenge, we develop a novel error decomposition that recasts decentralized GDRO as (i) a decentralized online convex optimization (D-OCO) problem and a prediction-with-expert-advice problem in SemiDec, and as (ii) two coupled D-OCO problems in Dec, with discrepancy terms capturing inconsistencies among local decisions.
By isolating consensus errors in analysis, we develop modular frameworks that support different consensus routines.
We then establish the \emph{first} per-worker convergence guarantees for optimizing the original GDRO objective over decentralized networks.
Instantiations with accelerated gossip and Choco-gossip achieve an $O(\sqrt{(\log m)/T})$ convergence rate over $T$ update rounds, while experiments demonstrate improved worst-worker performance over baselines.
\end{abstract}

\section{Introduction}
Decentralized optimization is a powerful framework for large-scale model training and has attracted considerable attention due to its network robustness, communication efficiency, and scalability~\citep{2018_COLA,2023_DFL_survey}.
Within this framework, $m$ workers exchange information over a communication network to collaboratively optimize a shared learning objective.
Classical methods train models by minimizing the average risk~\citep{2012_ADMM,Duchi12TAC,2016_DSGD}.
However, in practical deployments, workers' local data distributions are often highly heterogeneous~\citep{2023_survey_HFL}, which can lead to performance degradation across workers~\citep{pmlr-v119-hsieh20a}. As a result, minimizing only the average risk may be unfair to certain workers and lacks robustness under distribution shift~\citep{2019AFL}.

A promising solution is group distributionally robust optimization (GDRO)~\citep{Gouop_DRO_ICLR_2020}, which minimizes the maximum risk across workers.
For each worker $i\in[m]$, let $\mathcal{P}_i$ denote its local data distribution and define the local risk function as $R_i(\w)=\E_{\z\sim\mathcal{P}_i}[\ell(\w;\z)]$, where $\z$ denotes a random sample and $\ell(\cdot;\cdot)$ is the loss function.
GDRO solves the stochastic minimax problem
\begin{equation}\label{eq:gdro_intro:saddle}
\min_{\w\in\mathcal{W}} \max_{\q\in\Delta_m}
\left\{
\phi(\w,\q)=\sum_{i=1}^m q_i R_i(\w)
\right\},
\end{equation}
where $\Delta_m = \{\q\in\mathbb{R}^m_+ \mid \mathbf{1}^\top \q = 1\}$ denotes the simplex.
We focus on the convex setting and refer to $\w$ as the primal model and $\q$ as the dual vector.
Although GDRO has been extensively studied in centralized optimization~\citep{Online:Multiple:Distribution,zhang2023-SA-GDRO} and in server-based federated learning~\citep{2019AFL,2020_DRFA}, optimizing~(\ref{eq:gdro_intro:saddle}) over \emph{decentralized} networks remains challenging.
Existing decentralized GDRO methods only optimize a \emph{dual-regularized} variant of~(\ref{eq:gdro_intro:saddle}) with strongly concave regularization for the dual vector~\citep{2022_GDRO_dec,2023_GDRO_dec_cp}.
Consequently, their guarantees do not imply optimality for the worst-worker objective.
Moreover, these methods provide guarantees only for a globally averaged model, which requires additional communication to approximate.
Since deployment relies on locally maintained models, per-worker guarantees with respect to the global objective are more practically relevant~\citep{Duchi12TAC}.

In this paper, we directly optimize problem~(\ref{eq:gdro_intro:saddle}) over two representative decentralized architectures: (i) a \emph{semi-decentralized} (SemiDec) setting with a lightweight server and worker-to-worker communication, and (ii) a \emph{fully decentralized} (Dec) setting with only worker-to-worker communication.
Unlike decentralized optimization with a single variable~\citep{TAC_09_Dec}, GDRO couples the primal model and the dual vector, making its decentralized dynamics difficult to analyze.
A natural attempt is to adapt the centralized GDRO decomposition~\citep{zhang2023-SA-GDRO,Somav2}, which bounds the optimization error by the sum of primal and dual regret terms.
However, applying this decomposition to each worker's local model yields worker-dependent losses whose evaluation requires cross-worker information that is not locally accessible.
To overcome this challenge, we develop a novel error decomposition that rewrites each worker's local objectives as shared global objectives plus discrepancy terms arising from inconsistent local decisions.
This reformulation turns decentralized GDRO into tractable online learning subproblems.
In SemiDec, where a lightweight server maintains a shared dual vector, it yields a decentralized online convex optimization (D-OCO)~\citep{2013_DOCO_TKDE} problem for the primal updates and a prediction with expert advice (PEA)~\citep{Hedge_1997} problem for the dual updates, with model-discrepancy terms among workers.
In Dec, where each worker maintains both a local model and a local dual vector, it yields two coupled D-OCO problems together with primal and dual discrepancy terms.

Building on this decomposition, we further isolate communication effects through consensus-error terms in the regret and discrepancy bounds.
This separation decouples the optimization analysis from the specific consensus mechanism, yielding two modular frameworks that support different consensus routines.
Algorithmically, both employ follow-the-regularized-leader (FTRL)~\citep{orabona2023} for the primal and dual updates, while workers track the cumulative global gradient of each decentralized variable via the selected consensus routine.
We then establish convergence guarantees for each worker's time-averaged local model under the original GDRO objective.
By instantiating the frameworks with accelerated gossip (AG)~\citep{2011_liu_AG} and Choco-gossip (CG)~\citep{ICML19Cho_Gossip}, we obtain concrete algorithms with an $O(\sqrt{(\log m)/T})$ convergence rate over $T$ update rounds. Our main contributions are summarized as follows:
\begin{compactitem}
	\item We propose a novel error decomposition that reduces decentralized GDRO to a D-OCO problem and a PEA problem in SemiDec, and to two coupled D-OCO problems in Dec, with discrepancy terms capturing inconsistencies among local decisions.
    \item Building on the decomposition, we develop modular frameworks for GDRO in the SemiDec and Dec settings that accommodate different consensus routines, and establish per-worker convergence guarantees in terms of the resulting consensus errors.
	\item By instantiating the consensus routine with AG and CG, we obtain concrete algorithms with a per-worker convergence rate of $O(\sqrt{(\log m)/T})$, and empirically show improved worst-worker performance over prior methods.
\end{compactitem}

\section{Related Work}\label{app:related_work}
\textbf{Decentralized minimax optimization.}
Decentralized learning has been studied in SemiDec~\citep{2022_semi_Dec_HLSGD,2025_semi_Dec_PISCO} and Dec architectures~\citep{2007_DL_Saber,2012_DL_Chen}, with applications to energy systems~\citep{2013_Energy,2021_Ener_Semi}, wireless sensor networks~\citep{wireless_DL}, and autonomous vehicles~\citep{2024_Auto_Semi}.
Meanwhile, many machine learning problems can be formulated as minimax optimization~\citep{2017_RARL,2020_GNN}, motivating its study over decentralized networks. Most existing methods focus on regimes different from the convex-concave setting, including nonconvex-nonconcave~\citep{2020_NC_NC}, nonconvex-strongly-concave~\citep{2020_NC_SC,2021_NC_SC2,2025_DMMO}, and strongly-convex-strongly-concave objectives~\citep{2023_SC_SC}. 
Although \citet{2015_DL_OCO} study convex--concave problems, their dual variables are Lagrange multipliers for consensus constraints rather than distributional weights.

\textbf{Group distributionally robust optimization.}
GDRO~\citep{Gouop_DRO_ICLR_2020} minimizes the worst-group risk to learn models that perform uniformly well across multiple data distributions.
It has been studied in centralized settings~\citep{Gouop_DRO_ICLR_2020,Online:Multiple:Distribution,zhang2023-SA-GDRO,Bai2025FSQ,GDRO_TPAMI_2026}, server-based federated settings~\citep{2019AFL,2020_DRFA,2024_NIPS_FGDRO_Yang,2026_AFLCom_Bai}, and decentralized settings~\citep{2022_GDRO_dec,2023_GDRO_dec_cp} with a strongly concave regularizer.
Specifically, DR-DSGD~\citep{2022_GDRO_dec} adds Kullback--Leibler regularization to the dual vector, which transforms the regularized minimax problem into a single-level minimization problem, thus simplifying the algorithm design and analysis.
AD-GDA~\citep{2023_GDRO_dec_cp} further considers general strongly concave regularizers for the dual vector and solves the resulting problem using decentralized gradient descent-ascent with compressed consensus. 
Despite these advances, the original GDRO objective remains unexplored in decentralized settings.

\textbf{Decentralized online convex optimization.}
D-OCO~\citep{2013_DOCO_TKDE} considers multiple workers that collaboratively optimize a global objective in an online manner over a decentralized network, and algorithms typically include two key steps: decision updates and consensus. Existing update methods include decentralized online gradient descent~\citep{2013_DOCO_TKDE}, decentralized follow-the-regularized-leader (FTRL)~\citep{13_D_FTRL}, and their variants~\citep{2017_PF_DOCO,wan2022PFDOL,Wang_Wan_Zhang_Zhang_2023}. 
For information aggregation, most approaches adopt gossip for neighborhood averaging~\citep{XIAO2004_gossip}. 
\citet{Wan25JMLR_optimal_DOCO} introduce online accelerated gossip by combining AG~\citep{2011_liu_AG} with a blocking update mechanism, achieving nearly optimal regret bounds.
To improve communication efficiency, \citet{tu2022distributed} combine D-OGD with the difference-compressed consensus mechanism~\citep{ICML19Cho_Gossip} and \citet{2026_yang_DOCO} improve the results via a two-level blocking update framework. Subsequently, \citet{2026_wan_DOCO_COM} develop FTRL-type algorithms with compressed communication, matching the best known full-information regret bounds and improving both regret and communication costs in the bandit setting.

\section{Preliminaries} \label{sec:preliminaries}
We equip $\W$ with a closed, $1$-strongly convex function $\nu_w:\W\to\mathbb{R}$ with respect to $\|\cdot\|_2$, and $\Delta_m$ with negative entropy $\nu_q(\q)=\sum_{i=1}^m q_i\ln q_i$.
For any approximate solution $\w \in \W$ of problem~(\ref{eq:gdro_intro:saddle}), its optimization error is defined as $\varepsilon_{\phi}(\w)= \max_{i\in[m]} R_i(\w)-\min_{\w\in\W}\max_{i\in[m]} R_i(\w)$. We adopt standard assumptions as follows~\citep{2019AFL,zhang2023-SA-GDRO,2023_GDRO_dec_cp}.
\begin{assumption}\label{ass:domain_bound}
$\W$ is nonempty, compact, convex with $\max_{\w \in \W} \nu_w(\w) -\min_{\w \in \W} \nu_w(\w) \leq D^2$.
\end{assumption}
\begin{assumption}\label{ass:basic_pointwise}
For each $i\in[m]$, the risk $R_i(\cdot)$ is convex on $\W$, and for all $\w\in\W$, the stochastic gradient is unbiased, i.e., $\E_{\z\sim\P_i}[\nabla \ell(\w;\z)] = \nabla R_i(\w)$. In addition, for all $\w\in\W$ and all $\z$,
\begin{equation}\label{eqn:ass:value_gradient}
    0 \le \ell(\w;\z) \le 1,
    \qquad
    \|\nabla \ell(\w;\z)\|_2 \le G.
\end{equation}
\end{assumption}

\textbf{Communication network.}
The worker-to-worker communication network is represented by an undirected graph $\G=(\V,\mathcal{E})$, where $\V=[m]$ is the set of workers and $\mathcal{E}\subseteq \V\times\V$ denotes the set of communication links. 
An edge $(i, j) \in \mathcal{E}$ indicates that workers $i$ and $j$ can exchange information, and we assume $(i,i)\in\mathcal{E}$ for all $i\in\V$. 
Let $\mathcal{N}_i = \{ j \in \mathcal{V} : (i, j) \in \mathcal{E} \}$ denote the set of neighbors of worker $i$.
To model information mixing over $\G$, we associate it with a matrix $P\in\R^{m\times m}$, where $P_{i,j}>0$ if $(i,j)\in\mathcal{E}$ and $P_{i,j}=0$ otherwise. 
Following standard decentralized optimization~\citep{TAC_09_Dec,Duchi12TAC}, we impose the following assumption on $P$.
\begin{assumption} \label{ass:graph}
	The communication matrix $P$ satisfies $P \cdot \mathbf{1}_m = \mathbf{1}_m$, $P = P^\top$. Moreover, the graph $\mathcal{G}$ is connected, i.e., there exists a path between any two workers.
\end{assumption}
Under Assumption~\ref{ass:graph}, the largest eigenvalue of $P$ is $\lambda_1(P)=1$, while the remaining eigenvalues satisfy $\lambda_\star(P):=\max_{i\geq 2}|\lambda_i(P)|<1$~\citep{2012_matrix_analysis}.
We define the spectral gap as $\rho:=1-\lambda_\star(P)\in(0,1]$ and let $\beta:=\|I_m-P\|_2\in[0,2]$.

\textbf{Consensus routine.}\label{subsec:consensus}
We model the information-exchange mechanism as a black-box consensus routine $\A_c$. It returns updated states $\{\x^{(i)}_{L+1}, \hat{\x}_{L+1}^{(i)}\}_{i\in[m]}=\A_c(\{\x^{(i)}_{1}, \hat{\x}_{1}^{(i)}\}_{i\in[m]})$ after $L$ communication rounds, where $\x_{1}^{(i)}$ is the initial state of worker $i$, and $\hat{\x}_{t}^{(i)}$ denotes auxiliary states. 
We later instantiate $\A_c$ with two representative procedures, AG~\citep{2011_liu_AG} and CG~\citep{ICML19Cho_Gossip}, which respectively provide fast mixing and compressed communication. The full procedures and corresponding consensus-error bounds are provided in Appendix~\ref{app:consensus}.

\section{Semi-decentralized GDRO}\label{sec:semi_dec_GDRO}
We first consider SemiDec with a high-bandwidth worker-to-worker network and a low-bandwidth server--worker network.
Such architectures arise naturally in wireless systems, where nearby devices can exchange high-dimensional information through fast D2D links, while a base station provides global coordination over bandwidth-limited device-to-server links~\citep{2021_SemiDec_Lin,2022_semi_Dec_HLSGD}.
In SemiDec, we let the server only maintain the shared $m$-dimensional dual vector, while high-dimensional model update information is communicated over the worker-to-worker network.

\subsection{Error Decomposition} \label{sec:semi:error_dec}
Unlike prior work that provides guarantees only for a globally averaged model~\citep{2022_GDRO_dec,2023_GDRO_dec_cp}, we establish guarantees for each worker's \emph{local} model.
This is more relevant in decentralized settings, since obtaining the global average requires additional communication and is generally only approximate under finite communication.
Let $\w_t^{(i)}$ denote worker $i$'s model at round $t$, and define its output after $T$ iterations as $\wb_T^{(i)} = \frac{1}{T} \sum_{t=1}^T \w_t^{(i)}$. 
To decouple the effects of primal and dual updates, our starting point is the error decomposition used in centralized GDRO~\citep{zhang2023-SA-GDRO,Somav2}. In particular, the optimization error of each $\wb_T^{(i)}$ can be bounded as
\begin{equation}\label{eqn:proof:semi_dec:error_decom_informal}
	\begin{split}
		&\varepsilon_{\phi}(\wb_T^{(i)}) \\
		\leq  &\frac{1}{T} \underbrace{\left(\max_{\q\in \Delta_{m}}  \sum_{t=1}^T \phi(\w_t^{(i)},\q)- \sum_{t=1}^T \phi(\w_t^{(i)},\q_t)\right)}_{\mathcal{R}_q^{(i)}} + \frac{1}{T} \underbrace{\left(\sum_{t=1}^T \phi(\w_t^{(i)},\q_t) - \sum_{t=1}^T\phi(\w^*,\q_t)\right)}_{\mathcal{R}_w^{(i)}},
	\end{split}
\end{equation}
where $\w^*\in\argmin_{\w\in\W}\max_{i\in[m]}R_i(\w)$ is a global solution.
The term $\mathcal{R}_q^{(i)}$ is the regret of $\{\q_t\}_{t=1}^T$ against the loss sequence $\{\phi(\w_t^{(i)},\cdot)\}_{t=1}^T$ induced by worker $i$. 
Similarly, $\mathcal{R}_w^{(i)}$ is the regret of $\{\w_t^{(i)}\}_{t=1}^T$ measured on a shared objective sequence $\{\phi(\cdot,\q_t)\}_{t=1}^T$, with comparator $\w^*$.

\textbf{Analysis of $\mathcal{R}_q^{(i)}$.}
Since each worker $i\in[m]$ induces its own regret term $\mathcal{R}_q^{(i)}$, the server needs to select a \emph{single} dual sequence $\{\q_t\}_{t=1}^T$ that controls $\mathcal{R}_q^{(i)}$ \emph{uniformly} for all $i\in[m]$.
Constructing such a sequence directly, however, would require access to $R_j(\w_t^{(i)})$ for every $i,j\in[m]$.
For $j\neq i$, this risk can only be evaluated at worker $j$, who does not have access to worker $i$'s model $\w_t^{(i)}$.
To overcome this information barrier, we decompose each worker-dependent objective into a shared global term and a discrepancy term. Specifically, we rewrite
\begin{equation}\label{eqn:proof:semi_dec:error_decom_informal_1}
    \phi(\w_t^{(i)},\q)=\sum_{j=1}^m q_j R_j(\w_t^{(j)})+\sum_{j=1}^m q_j\bigl(R_j(\w_t^{(i)})-R_j(\w_t^{(j)})\bigr).
\end{equation}
The first term is independent of $i$ and therefore can serve as the global objective for updating the shared dual vector. Under Assumption~\ref{ass:basic_pointwise}, the second term can be bounded by $G\max_{j\in[m]}\|\w_t^{(i)}-\w_t^{(j)}\|_2$, which quantifies the maximum model disagreement between worker~$i$ and the rest of the network at round~$t$.
Consequently, each $\mathcal{R}_q^{(i)}$ can be bounded by the shared dual regret
\begin{equation} \label{eqn:proof:semi_dec:error_decom_informal_regq}
    \max_{\q\in \Delta_{m}}\sum_{t=1}^T\sum_{j=1}^m q_{j}R_j(\w_t^{(j)}) -\sum_{t=1}^T\sum_{j=1}^m q_{t,j}R_j(\w_t^{(j)}),
\end{equation}
plus a model-discrepancy term.
The shared regret in~(\ref{eqn:proof:semi_dec:error_decom_informal_regq}) leads to a PEA formulation at the server, where workers act as experts and the server updates $\q_t$ using stochastic estimates of their local risks.

\textbf{Analysis of $\mathcal{R}_w^{(i)}$.}
In contrast to $\mathcal{R}_q^{(i)}$, which requires further decomposition,
each primal regret term $\mathcal{R}_w^{(i)}$ in~(\ref{eqn:proof:semi_dec:error_decom_informal})
is already evaluated on the same objective sequence
$\{\phi(\cdot,\q_t)\}_{t=1}^T$, i.e.,
\begin{equation} \label{eqn:DOCO}
    \mathcal{R}_w^{(i)}=\sum_{t=1}^T 
    \left(\sum_{j=1}^m q_{t,j} R_j(\w_t^{(i)})  - \sum_{j=1}^m q_{t,j} R_j(\w^{*}) \right).
\end{equation}
To obtain uniform guarantees for all $\{\mathcal{R}_w^{(i)}\}_{i=1}^m$, we let workers jointly solve a D-OCO problem~\citep{2013_DOCO_TKDE}, with each worker computing its local gradient and using consensus to aggregate information across the network.

In summary, the above decomposition reduces the optimization error of each worker in SemiDec to a shared dual regret handled by PEA, a primal regret handled by D-OCO, and a model-discrepancy term.
We next develop the corresponding algorithmic framework.

\noindent \textbf{Remark.}
Different from standard uniform-weight D-OCO~\citep{2013_DOCO_TKDE}, our formulation involves time-varying dual weights.
A direct reduction would require rescaling and leads to loose bounds. Instead, we exploit $\q\in\Delta_m$ to avoid the extra dependence on $m$; see Appendix~\ref{sec:proof:semi_decGDRO:w}.

\subsection{Framework Design and Instantiation}
We first introduce the server-side update of the dual vector $\q_t$.
To control the shared dual regret~(\ref{eqn:proof:semi_dec:error_decom_informal_regq}), the server updates $\q_{t}$ to maximize the shared term $\sum_{j=1}^m q_{j} R_j(\w_t^{(j)})$, equivalently minimizing $\sum_{j=1}^m q_{j} (1-R_j(\w_t^{(j)}))$ under Assumption~\ref{ass:basic_pointwise}.
Accordingly, in each round $t$, each worker $i$ uploads only the scalar stochastic loss $\ell(\w^{(i)}_{t};\z_t^{(i)})$ to the server. 
The server then forms gradient estimator
\begin{equation}
	\g_{q,t} = [\g_{q,t}^{(1)}, \ldots, \g_{q,t}^{(m)}]^{\top}, \quad \text{with} \quad \g_{q,t}^{(i)}=1-\ell(\w^{(i)}_{t};\z_t^{(i)}) \label{eqn:stoch:grad:q:semi_dec}
\end{equation}
and computes $\q_{t+1}$ via FTRL~\citep{orabona2023} as
\begin{align}
	\q_{t+1}=\argmin_{\q\in \Delta_m} \left\{ \langle \eta_q \Lb_t, \q \rangle + \nu_q(\q)\right\}, \quad \text{with} \quad \Lb_t = \sum_{s=1}^t \g_{q,s} \label{eqn:update:qt} 
\end{align}
The server then broadcasts $q_{t+1,i}$ to worker $i$, who uses it to scale the local stochastic gradient.
\begin{algorithm}[t]
	\caption{\textsc{GDRO-SemiDec}}
	\begin{algorithmic}[1]
		\STATE Initialize $\w_1^{(i)}=\argmin_{\w \in \W} \nu_w(\w)$, $\Lb_0=\mathbf{0}_m$, $\q_1=\frac{1}{m}\mathbf{1}_m$, $\Fb^{(i)}_0=\mathbf{0}_d$, $\hFb_0^{(i)}= \emptyset$ (\textsc{AG}) \textbf{or} $\mathbf{0}_d$ (\textsc{CG}) $\forall i\in[m]$.
		\FOR{$t=1,2,\dots, T$} 
		\STATE \textbf{In parallel for all workers $i \in[m]$:}
		\STATE \quad Draw a sample $\z_t^{(i)}$ and send $\ell(\w_t^{(i)};\z_t^{(i)})$ to the server
		\STATE \quad Compute $\g_{w,t}^{(i)}=q_{t,i}\cdot\nabla\ell(\w_t^{(i)};\z_t^{(i)})$ 
		\STATE \quad $\{\Fb^{(i)}_{t},\hFb_t^{(i)}\}_{i\in[m]}= \A_c(\{\Fb^{(i)}_{t-1}+\g_{w,t}^{(i)},\hFb_{t-1}^{(i)}\}_{i\in[m]})$ \hfill \textit{// communication among workers}
		\STATE \quad Compute $\w_{t+1}^{(i)}$ via~(\ref{eqn:update:wt:dec})
		\STATE \textbf{On server:}
		\STATE \quad Construct $\g_{q,t}$ in (\ref{eqn:stoch:grad:q:semi_dec}) and update $\Lb_{t}=\Lb_{t-1}+ \g_{q,t}$
        \STATE \quad Compute $\q_{t+1}$ via~(\ref{eqn:update:qt}) and send $q_{t+1,i}$ to worker $i$ for all $i\in[m]$ 
		\ENDFOR
        \STATE \textbf{Return} $\wb_T^{(i)}=T^{-1}\sum_{t=1}^T\w_t^{(i)}$ at each worker $i\in[m]$
	\end{algorithmic} \label{alg:Semi_Dec_GDRO}
\end{algorithm}

To update the local models, each worker performs FTRL using an estimate of the \emph{global} cumulative gradient. 
Concretely, at round $t$, worker $i$ maintains $\Fb_{t-1}^{(i)}\in\mathbb{R}^d$ and auxiliary state variables $\hFb_{t-1}^{(j)}$ for $j\in\mathcal{N}_i$ when required by the consensus routine $\A_c$.
Each worker $i$ first computes the weighted stochastic gradient $\g_{w,t}^{(i)} = q_{t,i}\nabla \ell(\w_t^{(i)};\z_t^{(i)})$, and then applies $\mathcal{A}_c$ to the message $\Fb_{t-1}^{(i)}+\g_{w,t}^{(i)}$, yielding an updated estimate $\Fb_t^{(i)}$.
Finally, worker $i$ computes its local model via
\begin{equation}\label{eqn:update:wt:dec}
	\w_{t+1}^{(i)}	= \argmin_{\w \in \W} \left\{ \left\langle \eta_w \Fb_{t}^{(i)}, \w \right\rangle + \nu_w(\w) \right\}.
\end{equation}

Since consensus errors affect both the primal regret $\mathcal{R}_w^{(i)}$ and the model-discrepancy term arising from $\mathcal{R}_q^{(i)}$, we isolate its contribution in the corresponding bounds and treat the consensus routine $\mathcal{A}_c$ as a black box.
We present the resulting framework termed \textsc{GDRO-SemiDec} in Algorithm~\ref{alg:Semi_Dec_GDRO}, and establish its theoretical guarantees as follows.
\begin{theorem} \label{thm:Semi_Dec_GDRO}
	Under Assumptions~\ref{ass:domain_bound}--\ref{ass:graph}, set the step sizes as $\eta_w=\sqrt{2}mD/(\sqrt{T}G)$ and $\eta_q = \sqrt{2(\ln m)/T}$. 
    Suppose the consensus routine $\A_c$ ensures that $\sum_{i=1}^m\Fb_{t}^{(i)}=\sum_{i=1}^m\sum_{\tau=1}^{t}\g_{w,\tau}^{(i)}$ and $\E[\max_{i\in[m]}\|\Fb_{t}^{(i)}-\frac{1}{m}\sum_{r=1}^m\Fb_{t}^{(r)}\|_2]\leq\varepsilon_c$ for all $t\in[T]$. Then, Algorithm~\ref{alg:Semi_Dec_GDRO} guarantees
	\begin{equation*} \label{eqn:semi_dec:error:tem}
		\E\left[\varepsilon_{\phi}(\wb_T^{(i)})\right]
		\le
		4\sqrt{\frac{\ln m}{T}}+\sqrt{\frac{2}{T}}DG+\frac{7\sqrt{2}mD}{\sqrt{T}}\varepsilon_c, \qquad\forall i\in[m].
	\end{equation*}
\end{theorem}
\noindent\textbf{Remark.}
Theorem~\ref{thm:Semi_Dec_GDRO} establishes a convergence guarantee for each worker's time-averaged local model that only depends on the consensus error $\varepsilon_c$. 
This modularity allows $\mathcal{A}_c$ to be instantiated by different consensus routines, with different update rules or communication rounds per iteration.

We next instantiate $\mathcal{A}_c$ with AG and CG, and present the corresponding convergence guarantees. While the analysis is formulated through an online learning framework, the optimization problem itself is offline, and thus we allow multiple communication steps per iteration.

\begin{corollary}\label{cor:Semi_Dec_AG}
Under the conditions of Theorem~\ref{thm:Semi_Dec_GDRO}, let $\mathcal{A}_c$ be AG in Algorithm~\ref{alg:AG}, with $P\succeq0$.
With the communication-round choice $L=O((\log m)/\sqrt{\rho})$ specified in Appendix~\ref{app:proof:Semi_Dec_AG}, Algorithm~\ref{alg:Semi_Dec_GDRO} guarantees, for all $i\in[m]$, $\E[\varepsilon_{\phi}(\wb_T^{(i)})] \le (8\sqrt{2}DG+4)\sqrt{(\ln m)/T}.$
\end{corollary}
\begin{corollary}\label{cor:Semi_Dec_CG}
Under the conditions of Theorem~\ref{thm:Semi_Dec_GDRO}, let $\mathcal{A}_c$ be CG in Algorithm~\ref{alg:choco-gossip}.
With the communication-round choice $L=O((\log m)/(\delta\rho^2))$ specified in Appendix~\ref{app:proof:Semi_Dec_CG}, Algorithm~\ref{alg:Semi_Dec_GDRO} guarantees, for all $i\in[m]$, $\E[\varepsilon_{\phi}(\wb_T^{(i)})] \le (15\sqrt{2}DG+4)\sqrt{(\ln m)/T}.$
\end{corollary}

\section{Decentralized GDRO} \label{sec:dec_GDRO}
We now consider the Dec setting, where the lightweight server in SemiDec is removed and each worker $i$ maintains both a local dual vector $\q_t^{(i)}$ and local model $\w_t^{(i)}$.
This introduces dual inconsistency in addition to the model inconsistency already present in SemiDec.

\subsection{Error Decomposition}
Using again the centralized error decomposition, we bound the optimization error of $\wb_T^{(i)}$ as
\begin{equation*}\begin{split}
		&\varepsilon_{\phi}(\wb_T^{(i)}) \\
		\leq  &\frac{1}{T} \underbrace{\left(\max_{\q\in \Delta_{m}}  \sum_{t=1}^T \phi(\w_t^{(i)},\q)- \sum_{t=1}^T \phi(\w_t^{(i)},\q_t^{(i)})\right)}_{\mathcal{R}_q^{(i)}}
        + \frac{1}{T} \underbrace{\left(\sum_{t=1}^T \phi(\w_t^{(i)},\q_t^{(i)}) - \sum_{t=1}^T\phi(\w^*,\q_t^{(i)})\right)}_{\mathcal{R}_w^{(i)}}.
	\end{split}
\end{equation*}
The key difference to~(\ref{eqn:proof:semi_dec:error_decom_informal}) is that the dual sequence is no longer shared across workers.
In SemiDec, the shared dual sequence $\{\q_t\}_{t=1}^T$ makes $\mathcal R_w^{(i)}$ a regret on a shared objective sequence, so only $\mathcal R_q^{(i)}$ requires reformulation.
In Dec, however, the local dual sequence $\{\q_t^{(i)}\}_{t=1}^T$ also makes $\mathcal R_w^{(i)}$ worker-dependent.
Meanwhile, $\mathcal R_q^{(i)}$ faces the same cross-worker information barrier as in SemiDec, but must now be controlled for worker-specific dual sequences.
Accordingly, we reuse the SemiDec reformulation for $\mathcal R_q^{(i)}$ and derive an additional reformulation for $\mathcal R_w^{(i)}$.

\textbf{Analysis of $\mathcal{R}_q^{(i)}$.}
As in the SemiDec setting, we apply the decomposition in~(\ref{eqn:proof:semi_dec:error_decom_informal_1}) to $\phi(\w_t^{(i)},\q)$. This decomposes $\mathcal{R}_q^{(i)}$ into two parts. The first is a local dual regret term:
\begin{equation}\label{eqn:proof:dec:rq_regret_global}
\max_{\q\in\Delta_m}\sum_{t=1}^T \sum_{j=1}^m q_{j} R_j(\w_t^{(j)}) - \sum_{t=1}^T\sum_{j=1}^m q_{t,j}^{(i)} R_j(\w_t^{(j)}),
\end{equation}
which evaluates the local dual sequence $\{\q_t^{(i)}\}_{t=1}^T$ on the shared global objective $\sum_{j=1}^m q_j R_j(\w_t^{(j)})$, leading to a D-OCO formulation.
The second is a discrepancy term bounded by $2G\sum_{t=1}^T\max_{j\in[m]}\|\w_t^{(i)}-\w_t^{(j)}\|_2$, which measures cumulative disagreement among local models.

\textbf{Analysis of $\mathcal{R}_w^{(i)}$.}
This term introduces an additional difficulty specific to Dec.
Since $\q_t^{(i)}$ varies across workers, the objective
$\phi(\w,\q_t^{(i)})$ is no longer shared across workers.
To remove this worker dependence, we introduce the proxy vector $\hat{\q}_t=(q_{t,1}^{(1)},\ldots,q_{t,m}^{(m)})$ and rewrite
\begin{equation}\label{eqn:proof:dec:rw_decomp_informal}
\phi(\w,\q_t^{(i)})=\sum_{j=1}^m q_{t,j}^{(j)}R_j(\w)+\sum_{j=1}^m \bigl(q_{t,j}^{(i)}-q_{t,j}^{(j)}\bigr)R_j(\w).
\end{equation}
The first term on the right-hand side of~(\ref{eqn:proof:dec:rw_decomp_informal})
is shared across workers and induces the primal regret
\begin{equation}\label{eqn:proof:dec:rw_regret_global}
\sum_{t=1}^T \sum_{j=1}^m q_{t,j}^{(j)}R_j(\w_t^{(i)})
-\sum_{t=1}^T \sum_{j=1}^m q_{t,j}^{(j)}R_j(\w^*).
\end{equation}
We control this regret uniformly across workers by letting them jointly solve a D-OCO problem.

A key challenge is controlling the second term in~(\ref{eqn:proof:dec:rw_decomp_informal}), which measures inconsistency among local dual vectors and can be bounded by $\|\q_t^{(i)}-\hat{\q}_t\|_1$ under Assumption~\ref{ass:basic_pointwise}.
Unlike the shared dual vector in SemiDec, the proxy $\hat{\q}_t$ need not lie in $\Delta_m$, because its coordinates come from different workers' local dual vectors.
Consequently, the simplex-based argument used in SemiDec no longer applies, while the direct bound $\|\q_t^{(i)}-\hat{\q}_t\|_1\leq m$ incurs an undesirable factor of $m$.
To overcome this difficulty, we introduce a virtual simplex reference sequence $\{\tilde{\q}_t\}_{t\geq 1}$ driven by the averaged dual gradient and bound $\|\q_t^{(i)}-\hat{\q}_t\|_1\le\|\q_t^{(i)}-\tilde{\q}_t\|_1+\|\hat{\q}_t-\tilde{\q}_t\|_1.$ 
This construction allows both deviations to be controlled through the dual consensus error, avoiding the additional factor of $m$.

In summary, the optimization error of each worker's output is bounded by two D-OCO regret terms over shared global objectives, together with two discrepancy terms.

\subsection{Framework Design and Instantiation}
To control the two coupled D-OCO regrets, we design a consensus-based FTRL scheme in which each worker uses $\mathcal{A}_{c,w}$ and $\mathcal{A}_{c,q}$ to maintain and update local estimates $\Fb_t^{(i)}$ and $\Lb_t^{(i)}$ of the cumulative global gradients for the primal and dual updates, respectively.

We first introduce the update of the local dual vectors.
To optimize the shared dual objective in~(\ref{eqn:proof:dec:rq_regret_global}), we equivalently minimize
$\sum_{j=1}^m q_j(1-R_j(\w_t^{(j)}))$
under Assumption~\ref{ass:basic_pointwise}.
At round $t$, each worker $i$ forms a sparse stochastic contribution to the global dual gradient $\g_{q,t}^{(i)}=\e_i\bigl(1-\ell(\w_t^{(i)};\z_t^{(i)})\bigr),$ where $\e_i\in\R^m$ denotes the $i$-th standard basis vector.
The cumulative messages $\Lb_{t-1}^{(i)}+\g_{q,t}^{(i)}$ are then aggregated through $\mathcal{A}_{c,q}$ to update $\Lb_t^{(i)}$.
Finally, each worker $i$ performs the FTRL update
\begin{equation}\label{eqn:update:q:dec}
    \q_{t+1}^{(i)}
    =\argmin_{\q\in \Delta_m}\Bigl\{\bigl\langle \eta_{q}\Lb_t^{(i)}, \q \bigr\rangle+\nu_q(\q)\Bigr\}.
\end{equation}

The model update follows the same consensus-based FTRL scheme as in SemiDec, but now uses locally maintained dual weights.
To optimize the shared primal objective $\sum_{j=1}^m q_{t,j}^{(j)}R_j(\w)$ in~(\ref{eqn:proof:dec:rw_regret_global}), each worker $i$ forms the weighted stochastic gradient $\g_{w,t}^{(i)}=q_{t,i}^{(i)}\nabla \ell(\w_t^{(i)};\z_t^{(i)})$. The workers then aggregate the cumulative messages $\Fb_{t-1}^{(i)}+\g_{w,t}^{(i)}$ through $\mathcal{A}_{c,w}$, yielding updated estimates $\Fb_t^{(i)}$. Finally, each worker $i$ updates its local model via the FTRL step in~(\ref{eqn:update:wt:dec}).

We summarize the resulting framework, termed \textsc{GDRO-Dec}, in Algorithm~\ref{alg:Dec_GDRO} and establish its convergence guarantee in the following theorem.
\begin{algorithm}[t]
	\caption{\textsc{GDRO-Dec}}
	\begin{algorithmic}[1]
		\STATE Initialize $\w_1^{(i)}=\argmin_{\w \in \W} \nu_w(\w)$, $\Fb^{(i)}_0=\mathbf{0}_d$, $\Lb^{(i)}_0=\mathbf{0}_m$, $\q^{(i)}_1=\frac{1}{m}\mathbf{1}_m$, $\hFb_0^{(i)}=\hLb^{(i)}_0= \emptyset$ (\textsc{AG}) \textbf{or} $\hFb_0^{(i)}=\mathbf{0}_d, \hLb^{(i)}_0= \mathbf{0}_m$ (\textsc{CG}), $\forall i\in[m]$.
		\FOR{$t=1,2,\dots, T$} 
		\STATE \textbf{In parallel for all workers $i \in[m]$:}
		\STATE \quad Draw a sample $\z_t^{(i)}$ and compute $\g_{w,t}^{(i)}=q_{t,i}^{(i)}\nabla \ell(\w_t^{(i)};\z_t^{(i)}), \g_{q,t}^{(i)}=\e_i(1-\ell(\w_t^{(i)};\z_t^{(i)}))$
		\STATE \quad $\{\Fb^{(i)}_{t}, \hFb^{(i)}_t\}_{i\in[m]} = \A_{c,w}(\{\Fb^{(i)}_{t-1} + \g_{w,t}^{(i)}, \hFb^{(i)}_{t-1}\}_{i\in[m]})$ \hfill \textit{// primal consensus}
		\STATE \quad $\{\Lb^{(i)}_{t}, \hLb^{(i)}_t\}_{i\in[m]} = \A_{c,q}(\{\Lb^{(i)}_{t-1} + \g_{q,t}^{(i)}, \hLb^{(i)}_{t-1}\}_{i\in[m]})$ \hfill \textit{// dual consensus}
		\STATE \quad Update $\q_{t+1}^{(i)}$ in~(\ref{eqn:update:q:dec}) and  $\w_{t+1}^{(i)}$ in~(\ref{eqn:update:wt:dec})
		\ENDFOR
        \STATE \textbf{Return} $\wb_T^{(i)}=T^{-1}\sum_{t=1}^T\w_t^{(i)}$ at each worker $i\in[m]$
	\end{algorithmic}\label{alg:Dec_GDRO}
\end{algorithm}
\begin{theorem}\label{thm:Dec_GDRO_modular}
Under Assumptions~\ref{ass:domain_bound}--\ref{ass:graph}, set $\eta_w=\sqrt{2}mD/(G\sqrt{T})$ and $\eta_q = m\sqrt{2(\ln m)/T}$.
Suppose the consensus routines satisfy $\sum_i\Fb_t^{(i)}=\sum_i\sum_{\tau=1}^t\g_{w,\tau}^{(i)}$ and $\sum_i\Lb_t^{(i)}=\sum_i\sum_{\tau=1}^t\g_{q,\tau}^{(i)}$.
Writing $\bFb_t=m^{-1}\sum_i\Fb_t^{(i)}$ and $\bLb_t=m^{-1}\sum_i\Lb_t^{(i)}$, assume $\bigl(\E[\max_i\|\Fb_t^{(i)}-\bFb_t\|_2^2]\bigr)^{1/2}\le\varepsilon_{c,w}$ and $\bigl(\E[\max_i\|\Lb_t^{(i)}-\bLb_t\|_2^2]\bigr)^{1/2}\le\varepsilon_{c,q}$ for all $t\in[T]$.
Then, Algorithm~\ref{alg:Dec_GDRO} guarantees
\begin{equation}\label{eqn:dec:error:modular}
\E\left[\varepsilon_{\phi}(\wb_T^{(i)})\right]
\le O\left(\frac{C_1+C_2+C_3+\sqrt{\log m}}{\sqrt{T}}\right), \qquad \forall i\in[m],
\end{equation}
where $C_1$, $C_2$, and $C_3$ are consensus-dependent terms defined in Appendix~\ref{app:proof:Dec_GDRO_modular}. Under exact consensus, $C_1=1$ and $C_2=C_3=0$.
\end{theorem}
Similar to Theorem~\ref{thm:Semi_Dec_GDRO}, Theorem~\ref{thm:Dec_GDRO_modular} isolates the effect of communication through the consensus errors $\varepsilon_{c,w}$ and $\varepsilon_{c,q}$, allowing $\mathcal{A}_{c,w}$ and $\mathcal{A}_{c,q}$ to be instantiated independently. 
For example, one may use AG for the low-dimensional dual updates and CG for the high-dimensional model updates. 
For simplicity, we next present results for setting both $\mathcal{A}_{c,w}$ and $\mathcal{A}_{c,q}$ to either AG or CG.
\begin{corollary}\label{cor:Dec_AG}
Under the conditions of Theorem~\ref{thm:Dec_GDRO_modular}, let $\mathcal{A}_{c,w}$ and $\mathcal{A}_{c,q}$ be AG in Algorithm~\ref{alg:AG}, with $P\succeq0$.
With the communication-round choices $L_w,L_q=O\left((\log m)/\sqrt{\rho}\right)$ specified in Appendix~\ref{app:proof:Dec_AG}, Algorithm~\ref{alg:Dec_GDRO} guarantees, for all $i\in[m]$, $\E[\varepsilon_{\phi}(\wb_T^{(i)})]\le  (32DG+9\sqrt{2})\sqrt{(\ln m)/T}.$
\end{corollary}
\begin{corollary}\label{cor:Dec_CG}
Under the conditions of Theorem~\ref{thm:Dec_GDRO_modular}, let $\mathcal{A}_{c,w}$ and $\mathcal{A}_{c,q}$ be CG in Algorithm~\ref{alg:choco-gossip}.
With the communication-round choices $L_w,L_q=O((\log m)/(\delta\rho^2))$ specified in Appendix~\ref{app:proof:Dec_CG}, Algorithm~\ref{alg:Dec_GDRO} guarantees, for all $i\in[m]$, $\E[\varepsilon_{\phi}(\wb_T^{(i)})]\le  (61DG+9\sqrt{2})\sqrt{(\ln m)/T}.$
\end{corollary}

\section{Discussion of Theoretical Guarantees}
\label{sec:theory_discussion}
To the best of our knowledge, Theorems~\ref{thm:Semi_Dec_GDRO} and~\ref{thm:Dec_GDRO_modular} provide the \emph{first} per-worker convergence guarantees for the original GDRO objective over decentralized networks. 
In both settings, the per-worker optimization error bound consists of an intrinsic optimization term and communication-induced consensus errors, which we discuss below.

First, the optimization term is tight in sample complexity up to a logarithmic factor in $m$. Under exact consensus, i.e., $\varepsilon_c=0$ in SemiDec and $\varepsilon_{c,w}=\varepsilon_{c,q}=0$ in Dec, the error reduces to $O(\sqrt{(\log m)/T})$. Thus, finding an $\varepsilon$-optimal solution requires $O((\log m)/\varepsilon^2)$ updates and $O((m\log m)/\varepsilon^2)$ samples, matching the best-known centralized GDRO results~\citep{zhang2023-SA-GDRO} and the information-theoretic lower bound $\Omega(m/\varepsilon^2)$ up to a factor of $O(\log m)$~\citep{2022_Soma_DRO}.

Second, to preserve the $O(\sqrt{(\log m)/T})$ rate, we choose the number of communication rounds in Corollaries~\ref{cor:Semi_Dec_AG}--\ref{cor:Dec_CG} according to the best-known scaling of the D-OCO subproblems. Specifically, the AG instantiations match the communication-round scaling of the near-optimal D-OCO method in~\citet{Wan25JMLR_optimal_DOCO}, while the CG instantiations inherit the state-of-the-art dependence on network connectivity and compression~\citep{2026_yang_DOCO}. Importantly, these choices are sufficient rather than necessary, and the main theorems remain valid for different numbers of consensus rounds.

\section{Experiments} \label{sec:experiments}
We evaluate our SemiDec and Dec methods on three datasets. Additional experimental details and worst-worker loss results are provided in Appendix~\ref{app:experiment}.

\noindent\textbf{Experimental setup.} We consider Adult~\citep{adult_dataset}, MNIST~\citep{MNIST}, and CIFAR-10~\citep{CIFAR10}.
We instantiate \textsc{GDRO-SemiDec} and \textsc{GDRO-Dec} with AG and CG, denoting the resulting variants by AG$(L)$ and CG$(L)$ within each setting, where $L$ is the number of consensus rounds per update. 
For Dec, we set $L_w=L_q=L$ and compare with DR-DSGD~\citep{2022_GDRO_dec} and AD-GDA~\citep{2023_GDRO_dec_cp}. 
AD-GDA and the CG-based variants use Top-$k$ compression~\citep{2018_SparSGD_toprandk} with a $50\%$ compression ratio.
Let $\mathrm{Acc}_j(\bar{\w}_t^{(i)})$ and $R_j(\bar{\w}_t^{(i)})$ denote the test accuracy and loss of worker $i$'s local model on worker $j$'s test data, respectively. We report worst-worker accuracy $\min_{i,j\in[m]}\mathrm{Acc}_j(\bar{\w}_t^{(i)})$ in the main text and worst-worker loss $\max_{i,j\in[m]}R_j(\bar{\w}_t^{(i)})$ in Appendix~\ref{app:exp_risk}.
Workers communicate over an undirected ring graph. We train linear classifiers on Adult and MNIST. On CIFAR-10, all methods use the same pretrained ResNet-20~\citep{polivin2025resnet20} with a frozen backbone and a randomly initialized linear head.
\begin{figure*}[t]
    \centering
    \includegraphics[width=\textwidth]{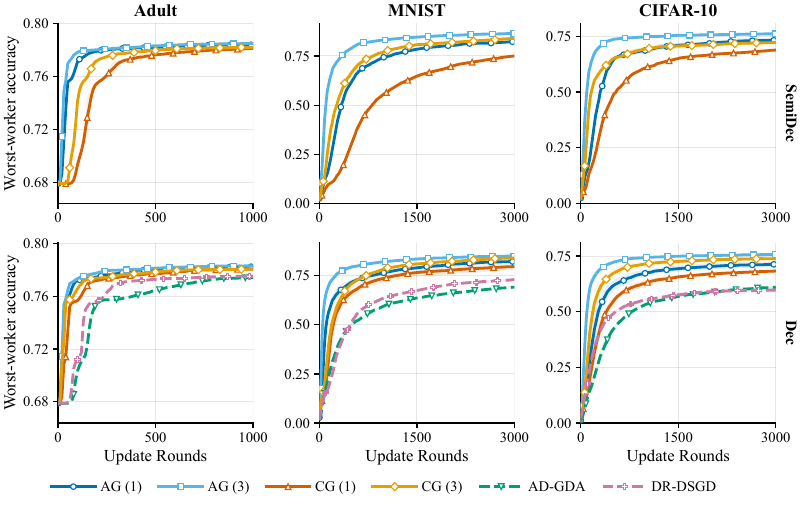}
    \caption{Worst-worker accuracy over update rounds on Adult, MNIST, and CIFAR-10. The top and bottom rows correspond to SemiDec and Dec, respectively. Curves report the mean over five random seeds, with shaded regions showing one standard deviation.}
    \label{fig:main_results}
\end{figure*}

\noindent\textbf{Results.}
Figure~\ref{fig:main_results} reports the worst-worker accuracy over update rounds in both settings.
First, in Dec, our variants generally converge faster and attain higher accuracy than the baselines.
These results support the practical advantage of directly optimizing the unregularized GDRO objective, whereas the baselines optimize dual-regularized surrogates.
Second, across both settings, the $L=3$ variants generally improve faster than their $L=1$ counterparts. This behavior is consistent with the consensus-error dependence in Theorems~\ref{thm:Semi_Dec_GDRO} and~\ref{thm:Dec_GDRO_modular}, as additional consensus rounds reduce the consensus errors.
Finally, both AG and CG exhibit stable convergence across the datasets, demonstrating the flexibility of our framework in accommodating different consensus routines.

\section{Conclusion}
In this paper, we study GDRO over decentralized networks in both SemiDec and Dec settings. 
Our key technical contribution is an error decomposition that bounds each worker's local optimization error by global regret terms plus disagreement-induced discrepancies, leading to a D-OCO problem and a PEA problem in SemiDec, and two coupled D-OCO problems in Dec.
Based on this decomposition, we develop modular frameworks supporting different consensus routines and bound the expected optimization error of each worker's time-averaged local model for the original GDRO objective. 
We then instantiate the frameworks with AG and CG, yielding an $O(\sqrt{(\log m)/T})$ convergence rate. Experiments demonstrate improved worst-worker performance over baselines.

\bibliography{ref}
\bibliographystyle{iclr2027_conference}

\appendix
\section{Consensus Routine Instantiations} \label{app:consensus}
In this appendix, we present two representative instantiations of the black-box consensus routine $\A_c$: accelerated gossip (AG) and Choco-gossip (CG). We provide their procedures and the standard contraction guarantees used in our analysis.

\subsection{Accelerated Gossip} \label{app:ag}
AG~\citep{2011_liu_AG} is an accelerated linear consensus scheme. Starting from $\x_0^{(i)}=\x_1^{(i)}$ for all $i\in[m]$, it performs $L$ communication rounds. At each round $k\in[L]$, worker $i$ exchanges $\x_k^{(i)}$ with its neighbors and updates
\[
\x_{k+1}^{(i)}=(1+\theta)\sum_{j\in\mathcal{N}_i}P_{i,j}\x_k^{(j)}-\theta \x_{k-1}^{(i)}.
\]
where $\theta>0$ is the stepsize. 
For the AG instantiation, we additionally assume that $P\succeq 0$. 
The AG procedure is summarized in Algorithm~\ref{alg:AG} and enjoys the following guarantee.
\begin{lemma}~\citep{2023_AG_the} \label{lem:ag_convergence}
Let the step size $\theta=\frac{1}{1+\sqrt{1-\lambda_\star^2(P)}}$.
Define $X_k= \big[\x_k^{(1)},\ldots,\x_k^{(m)}\big]\in\mathbb{R}^{d\times m}$ and the average $\bar{\x}= \frac{1}{m} X_0\mathbf{1}\in\mathbb{R}^d$. Then $\bar{\x}=\frac{1}{m} X_{L+1}\mathbf{1}$ and
\[
\bigl\|X_{L+1}-\bar{\x}\mathbf{1}^{\top}\bigr\|_F
\le
\sqrt{14}c^{L}
\bigl\|X_{0}-\bar{\x}\mathbf{1}^{\top}\bigr\|_F,
\]
where $c=1-(1-1/\sqrt{2})\sqrt{1-\lambda_\star(P)}$ and $\lambda_\star(P)$ denotes the second largest eigenvalue of $P$.
\end{lemma}
\begin{algorithm}[t]
\caption{AG($\{\x_1^{(i)}\}_{i\in[m]}$): Accelerated Gossip}
\label{alg:AG}
\textbf{Input:} Initial states $\{\x_1^{(i)}\}_{i\in[m]}$, rounds $L$, communication matrix $P\in\mathbb{R}^{m\times m}$, mixing coefficient $\theta>0$.

\begin{algorithmic}[1]
\STATE Initialize $\x_{0}^{(i)} = \x_{1}^{(i)}$ for all $i\in[m]$
\FOR{$k=1,\dots,L$}
  \STATE \textbf{In parallel for all workers $i\in[m]$:}
  \STATE Send $\x_k^{(i)}$ to all neighbors $j\in\mathcal{N}_i$ and receive $\{\x_k^{(j)}\}_{j\in\mathcal{N}_i}$
  \STATE Update $\x_{k+1}^{(i)} = (1+\theta) \sum_{j\in \mathcal{N}_i} P_{i,j}\ \x_k^{(j)} - \theta \x_{k-1}^{(i)}$
\ENDFOR
\STATE \textbf{Return} $\{\x_{L+1}^{(i)}\}_{i\in[m]}$
\end{algorithmic}
\end{algorithm}

\subsection{Choco-gossip} \label{app:cg}
CG is a communication-efficient consensus procedure that reduces bandwidth by exchanging \emph{compressed} messages. It relies on a compressor $\mathcal{C}$ defined as follows.
\begin{definition} \label{def:compressor}
Let $\mathcal{Q}\subseteq\mathbb{R}^d$ and $\delta\in(0,1]$.
An operator $\mathcal{C}:\mathbb{R}^d\to\mathbb{R}^d$ is a $\delta$-approximate compressor over $\mathcal{Q}$ if
\begin{equation}\label{eqn:define:com:sto}
\E\left[\left\|\mathcal{C}(\x)-\x\right\|_2^2\right]\leq(1-\delta)\|\x\|_2^2,
\qquad \forall\x\in\mathcal{Q}.
\end{equation}
where the expectation is taken with respect to the internal randomness of $\mathcal{C}$.
\end{definition}
At each communication round $k\in[L]$, CG maintains both the current states $\{\x_k^{(i)}\}_{i\in[m]}$ and the transmitted/received information $\{\hat{\x}_k^{(i)}\}_{i\in[m]}$ (i.e., worker $i$ stores $\{\hat{\x}_k^{(j)}\}_{j\in\mathcal{N}_i}$). Each worker $i$ sends the compressed difference $\Delta_k^{(i)}=\mathcal{C}(\x_k^{(i)}-\hat{\x}_k^{(i)})$ to its neighbors.
Upon receiving $\{\Delta_k^{(j)}\}_{j\in\mathcal{N}_i}$, worker $i$ updates $\{\hat{\x}_{k+1}^{(j)}\}_{j\in\mathcal{N}_i}$ and $\x_{k+1}^{(i)}$.
CG is summarized in Algorithm~\ref{alg:choco-gossip} and satisfies the following guarantee.
\begin{lemma} \citep{ICML19Cho_Gossip} \label{lem:choco-gossip}
	Let $P$ satisfy Assumption~\ref{ass:graph}, and assume $\mathcal{C}$ satisfies~(\ref{eqn:define:com:sto}) on $\mathbb{R}^d$ with fresh internal randomness at each call. Define  
	$X_t = [\x_t^{(1)}, \ldots, \x_t^{(m)}], \hat{X}_t = [\hat{\x}_t^{(1)}, \ldots, \hat{\x}_t^{(m)}],$	and $\bar{\x}_t = \frac{1}{m} \sum_{i=1}^m \x_t^{(i)}$. 
	Algorithm~\ref{alg:choco-gossip} ensures \( \bar{\x}_t = \bar{\x} \) for all \( t\in[L] \). Define \(\bar{X}= \bar{\x}\mathbf{1}^\top\), and 
	$e_{t+1} = \E [ \| X_{t+1} - \hat{X}_{t+1} \|_F^2 + \| X_{t+1} - \bar{X} \|_F^2 ]$.
	With step size $\gamma = \frac{\rho\delta}{16\rho+\rho^{2}+4\beta^{2}+2\rho\beta^{2}-8\rho\delta}$, it holds that
	\begin{equation} \label{eqn:lem:ICLR25}
		e_{L+1} \leq \left(1 - \frac{\gamma \rho}{2} \right) e_L \leq \left(1 - \frac{\gamma \rho}{2} \right)^L e_1.
	\end{equation}
\end{lemma}
\begin{algorithm}[t]
	\caption{CG($\{\x_{1}^{(i)}, \hat{\x}_{1}^{(i)}\}_{i\in[m]}$): Choco-gossip}
	\textbf{Input:} States $\x_{1}^{(i)}, \hat{\x}_{1}^{(i)}$ for all $i \in [m]$, step size $\gamma$, rounds $L$, compressor $\C$, communication matrix $P$
    
	\begin{algorithmic}[1]
		\FOR{$k=1,2,\dots, L$}
		\STATE \textbf{In parallel for all workers $i \in[m]$:}
		\STATE Send $\Delta^{(i)}_{k}=\C (\x_{k}^{(i)}-\hat{\x}_{k}^{(i)})$ to all neighbors $j\in\mathcal{N}_i$
		\STATE Update $\hat{\x}_{k+1}^{(j)}=\hat{\x}_{k}^{(j)}+\Delta^{(j)}_{k}$ for $j\in\mathcal{N}_i$
		\STATE Update $\x_{k+1}^{(i)} = \x_{k}^{(i)} + \gamma\sum_{j\in \mathcal{N}_i}P_{i,j} (\hat{\x}_{k+1}^{(j)}-\hat{\x}_{k+1}^{(i)})$ 
		\ENDFOR
		\STATE \textbf{Return} $\{\x_{L+1}^{(i)},\hat{\x}_{L+1}^{(i)}\}_{i\in[m]}$
	\end{algorithmic} \label{alg:choco-gossip}
\end{algorithm}

\section{More Experiment Results}\label{app:experiment}
\subsection{Experimental Details} \label{app:experiment_details}
We run all experiments in Python on an Ubuntu 20.04.4 LTS server with two Intel(R) Xeon(R) Platinum 8358P processors at 2.60\,GHz.

For Adult, we use the official training and test files and form $m=6$ workers using the Cartesian product of race $\in\{\mathrm{White},\mathrm{Black},\mathrm{Other}\}$ and sex $\in\{\mathrm{Male},\mathrm{Female}\}$. For MNIST and CIFAR-10, we use the official training and test sets with $m=10$ workers. Each split is partitioned independently: $80\%$ of samples from class $i$ are assigned to worker $i$, and the remaining $20\%$ are distributed approximately uniformly among the other workers.

Workers communicate over an undirected ring with mixing matrix $P=I-\frac{L_G}{2d_{\max}}$.
For the degree-two ring used in our experiments, $P_{ii}=1/2$ and $P_{i,i\pm1}=1/4$, with indices taken modulo $m$. The resulting matrix is symmetric, doubly stochastic, and positive semidefinite for both $m=6$ and $m=10$, and therefore satisfies the matrix condition in our AG corollaries.

We use a linear logistic classifier on Adult and a linear $784$-to-$10$ classifier on MNIST. On Adult, the trainable parameter vector is projected onto the Euclidean ball of radius $3$ after each update. For CIFAR-10, we use a pretrained ResNet-20 as a fixed feature extractor. All backbone parameters and BatchNorm statistics are frozen, and only a randomly initialized $64$-to-$10$ linear head is trained.

For the proposed methods, the primal update uses the original batch loss, while the dual feedback is clipped and normalized. Let $\ell_{\mathrm{ce}}(\w;z)$ denote the cross-entropy loss and define
\[
\ell_C(\w;z)=\frac{\min\{\ell_{\mathrm{ce}}(\w;z),C\}}{C},
\qquad C=5.
\]
The dual update uses
\[
\hat{s}_{t,i}=1-\ell_C(\w_t^{(i)};z_t^{(i)}),
\]
so that $\hat{s}_{t,i}\in[0,1]$. The primal update and evaluation use the original unclipped loss.

Top-$k$ compression retains the $\lceil0.5d\rceil$ coordinates with largest absolute magnitude in each parameter tensor. We use a mini-batch size of $128$ for Adult and MNIST and $32$ for CIFAR-10. For AD-GDA, we use its chi-square dual regularization with $\alpha=0.1$, $0.1$, and $0.05$ on Adult, MNIST, and CIFAR-10, respectively; for DR-DSGD, we set $\mu=3$, $3$, and $1$, respectively.

\subsection{Worst-Worker Loss}\label{app:exp_risk}
\begin{figure*}[t]
    \centering
    \includegraphics[width=0.9\textwidth]{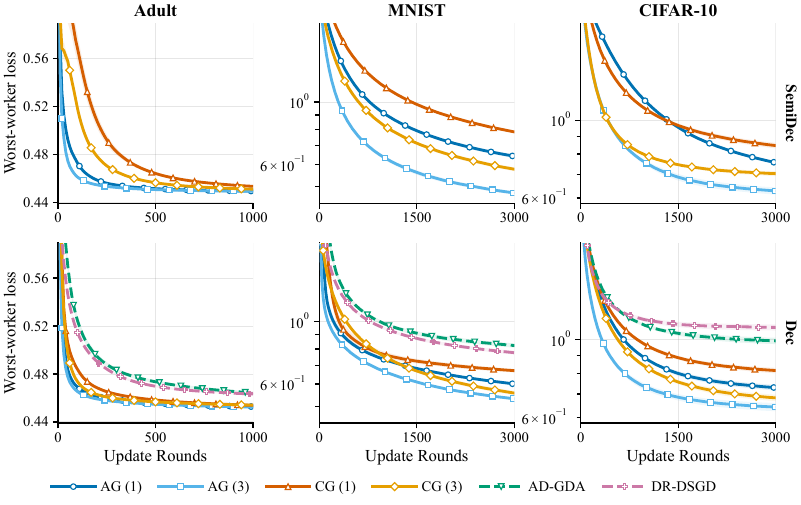}
    \caption{Worst-worker loss over update rounds on Adult, MNIST, and CIFAR-10. The top and bottom rows correspond to SemiDec and Dec, respectively. Parenthesized numbers denote the consensus rounds per update.}
    \label{fig:additional_risk}
\end{figure*}
Figure~\ref{fig:additional_risk} reports the loss curves corresponding to the accuracy results in Section~\ref{sec:experiments}. In SemiDec, worst-worker loss decreases for all variants. AG decreases more rapidly than CG on MNIST, while the $L=3$ variants generally improve earlier, particularly on CIFAR-10. In Dec, our variants generally finish below AD-GDA and DR-DSGD, with the clearest separation on CIFAR-10 and smaller differences on Adult and MNIST.

\section{Further Discussion of Theoretical Guarantees}
\label{app:theory_discussion}
Let $|\mathcal E_{\mathrm{ns}}|$ denote the number of non-self undirected edges and let $N$ count bidirectional neighbor exchanges. Over $T$ updates, SemiDec uses $N_{\mathrm{SemiDec}}=|\mathcal E_{\mathrm{ns}}|TL$ exchanges, whereas Dec uses $N_{\mathrm{Dec}}=|\mathcal E_{\mathrm{ns}}|T(L_w+L_q)$. For AG, $L,L_w,L_q=O((\log m)/\sqrt{\rho})$, giving $N_{\mathrm{AG}}=O(|\mathcal E_{\mathrm{ns}}|(\log m)^2/(\varepsilon^2\sqrt{\rho}))$ to reach expected error $\varepsilon$, up to a constant factor between SemiDec and Dec. For CG, $L,L_w,L_q=O((\log m)/(\delta\rho^2))$, yielding $N_{\mathrm{CG}}=O(|\mathcal E_{\mathrm{ns}}|(\log m)^2/(\varepsilon^2\delta\rho^2))$. 

Theoretically, the need for multiple consensus rounds is inherited from the underlying D-OCO subproblems. In particular, \citet{Wan25JMLR_optimal_DOCO} show that using a single AG step per round leads to a worse regret bound, whereas their blocking-based method uses $O((\log m)/\sqrt{\rho})$ communication rounds between consecutive decision updates and achieves a near-optimal regret bound, matching the corresponding lower bound up to a logarithmic factor. Under the correspondence between one online block and one update in our offline framework, this gives exactly the same $O((\log m)/\sqrt{\rho})$ communication-round order as our AG instantiations. The CG instantiations similarly inherit the state-of-the-art network and compression dependence of the corresponding compressed D-OCO method~\citep{2026_yang_DOCO}.

From a communication perspective, SemiDec and Dec require the same order of peer-to-peer consensus rounds for the $d$-dimensional primal information. Dec additionally communicates $m$-dimensional dual information through worker-to-worker consensus, whereas SemiDec maintains the shared dual vector through scalar server--worker communication.

\section{Analysis}   \label{app:analysis}
\subsection{Useful Lemmas}
We first state two lemmas for FTRL that will be used in our analysis. The first lemma controls the regret in terms of the norm of the gradients.
\begin{lemma}\citep[Corollary 7.9.]{orabona2023} \label{lem:FTRL_Orabona}
	Let $V$ be convex and let $\psi:V\to\mathbb{R}$ be a closed and $\mu$-strongly convex function w. r. t. $\| \cdot \|$. Set the sequence of regularizers as $\psi_t(\boldsymbol{x})=\frac1{\eta_{t-1}}(\psi(\boldsymbol{x})-\min_{\boldsymbol{z}\in V}\psi(\boldsymbol{z}))$, where $\eta_0\ge\eta_1\ge\cdots\ge\eta_T>0$ and $\boldsymbol{x}_t\in\operatorname{argmin}_{\boldsymbol{x}\in V}\psi_t(\boldsymbol{x})+\sum_{i=1}^{t-1}\ell_i(\boldsymbol{x})$. Assume $\ell_t$ is convex and closed, and let $\boldsymbol{g}_t\in\partial\ell_t(\boldsymbol{x}_t)$. Then, for every $\boldsymbol{u}\in V$, FTRL guarantees
	$$\sum_{t=1}^T\ell_t(\boldsymbol{x}_t)-\sum_{t=1}^T\ell_t(\boldsymbol{u})\leq\frac{\psi(\boldsymbol{u})-\min_{\boldsymbol{x}\in V}\psi(\boldsymbol{x})}{\eta_{T-1}}+\frac1{2\mu}\sum_{t=1}^T\eta_{t-1}\|\boldsymbol{g}_t\|_*^2.$$
\end{lemma}
The second lemma is a local-norm instantiation of Lemma~7.14 by \citet{orabona2023}, specialized to the simplex and the negative-entropy regularizer (see \S~7.5 in \citet{orabona2023} for further details).
\begin{lemma}\citep[Lemma 7.14]{orabona2023} \label{lem:FTRL_Orabona:local_norm}
	Denote by $F_t(\boldsymbol{x})=\psi_t(\boldsymbol{x})+\sum_{i=1}^{t-1}\left\langle\boldsymbol{g}_i, \boldsymbol{x}\right\rangle$ with $\psi_t(\boldsymbol{x})=\frac1{\eta_{t}}(\sum_{i=1}^mx_i\ln x_i-\min_{\boldsymbol{x}\in\Delta_m}\sum_{i=1}^mx_i\ln x_i)$, where $\eta_{t+1}\leq\eta_t,\forall t=1,\ldots,T$. Set $\boldsymbol{x}_t\in\operatorname*{argmin}_{\boldsymbol{x}\in \Delta_m}$ $F_t(\boldsymbol{x})$. Assume $\boldsymbol{g}_t\in\mathbb{R}_+^m$ and $\eta_t>0$ for all $t$. The following inequality holds for any $\boldsymbol{u} \in \Delta_m$,
	$$
	\sum_{t=1}^T\left\langle\boldsymbol{g}_t, \boldsymbol{x}_t-\boldsymbol{u}\right\rangle \leq \psi_{T}(\boldsymbol{u})-\min _{\boldsymbol{x} \in \Delta_m} \psi_1(\boldsymbol{x})+\frac{1}{2} \sum_{t=1}^T \eta_{t}\sum_{i=1}^m g_{t,i}^2x_{t,i}.
	$$
\end{lemma}

\begin{lemma}\citep[Lemma 5]{Duchi12TAC} \label{lem:Projection}
	Let $\mathcal{X}$ be nonempty, closed, and convex, let $\psi$ be closed and $1$-strongly convex with respect to $\|\cdot\|$, and let $\alpha>0$. Define $\Pi_{\mathcal{X}}^\psi(z,\alpha)=\underset{x\in\mathcal{X}}{\operatorname*{\mathrm{argmin}}}\left\{\langle z,x\rangle+\frac{1}{\alpha}\psi(x)\right\}$. For arbitrary $u, v \in \mathbb{R}^d$, we have
	\begin{equation}\label{eqn:lemma:projection}
		\left\|\Pi_{\mathcal{X}}^{\psi}(u, \alpha) - \Pi_{\mathcal{X}}^{\psi}(v, \alpha)\right\| \leq \alpha \|u - v\|_*.
	\end{equation}
\end{lemma}
\begin{lemma} \citep[Lemma A.7]{cesa2006prediction}  \label{lem:azuma}
	Let $V_1, V_2,  \ldots$  be a martingale difference sequence with respect to some sequence $X_1, X_2, \ldots$ such that $V_i \in [A_i , A_i + c_i ]$ for some random variable $A_i$, measurable with respect to $X_1, \ldots , X_{i-1}$ and a positive constant $c_i$. If $S_n = \sum_{i=1}^n V_i$, then for any
	$t > 0$,
	\[
	\Pr[ S_n > t] \leq \exp \left( -\frac{2t^2}{\sum_{i=1}^n c_i^2} \right).
	\]
\end{lemma}
\begin{lemma} \citep[\S~3.2]{Bandit:suvery} \label{lemma:high_pro_bound_to_Exp}
	For any real-valued random variable $X$, $\E[X] \leq \int_{0}^1 \frac{1}{\delta} \Pr \left[X > \ln \frac{1}{\delta}\right]~\mathrm{d} \delta.$
\end{lemma}

Based on Lemmas~\ref{lem:ag_convergence} and~\ref{lem:choco-gossip}, we obtain the following bound on the communication error.
\begin{lemma}[Accelerated Gossip Error Bound]\label{lem:iter_bound:ag_gossip}
Let $P\succeq0$ satisfy Assumption~\ref{ass:graph}, let $\rho=1-\lambda_\star(P)$, and initialize $\Xb_0^{(i)}=\mathbf{0}$ for all $i\in[m]$.
Consider the Accelerated Gossip operator
\begin{equation} \label{lem:acc_gossip_update}
\{\Xb^{(i)}_{t}\}_{i\in[m]}
= \text{\rm AG}\big(\{\Xb^{(i)}_{t-1} + \g_{t}^{(i)}\}_{i\in[m]}, L, P\big),
\end{equation}
i.e., run Algorithm~\ref{alg:AG} for $L$ communication rounds with $\theta= \frac{1}{1+\sqrt{1-\lambda_\star^2(P)}}$. 
Let $\Xb_t=[\Xb^{(1)}_{t},\ldots,\Xb^{(m)}_{t}]$ and
$\bXb_t = \frac{1}{m}\Xb_t \mathbf{1}$ be the column-wise average and denote $c=1-\Bigl(1-\tfrac{1}{\sqrt{2}}\Bigr)\sqrt{\rho}$.
Fix a constant $k>2$ and choose the number of gossip rounds as
\[
L = \left\lceil \frac{\ln(14k)}{2\left(1-\frac{1}{\sqrt{2}}\right)\sqrt{\rho}} \right\rceil,
\]
so that $14c^{2L}\le \frac{1}{k}$.
Assume additionally that $\sum_{i=1}^m \|\g_t^{(i)}\|_2^2 \le Q$ for all $t\in[T]$.
Then the per-round communication error satisfies
\begin{equation}\label{eq:lem:ag_gossip}
\|\Xb_t-\bXb_t\mathbf{1}^\top\|_F
\le \frac{\sqrt{Q}}{\sqrt{k}-1},\qquad \forall t\in[T].
\end{equation}
Moreover, the average is preserved in the sense that
$\bXb_t = \bXb_{t-1} + \bg_t$, where $\g_t=[\g_t^{(1)},\ldots,\g_t^{(m)}]$
and $\bg_t = \frac{1}{m}\g_t\mathbf{1}$.
\end{lemma}

\begin{lemma}[Choco-gossip Error Bound] \label{lem:semi_dec:choco-gossip}
	Let $\C$ satisfy~(\ref{eqn:define:com:sto}) on $\mathbb{R}^d$ with fresh internal randomness at each call, and let $P$ satisfy Assumption~\ref{ass:graph}. Initialize $\Xb_0^{(i)}=\hXb_0^{(i)}=\mathbf{0}$ for all $i\in[m]$. Consider the Choco-gossip operator 
    \begin{equation} \label{lem:choco_gossip_update}
        \{\Xb^{(i)}_{t},\hXb^{(i)}_t\}_{i\in[m]} = \text{\rm CG}(\{\Xb^{(i)}_{t-1} + \g_{t}^{(i)},\hXb^{(i)}_{t-1}\}_{i\in[m]}, \gamma, L, \C, P)    
    \end{equation}
    and let the gossip stepsize $\gamma = \frac{\delta \rho}{16\rho+\rho^{2}+4\beta^{2}+2\rho\beta^{2}-8\rho\delta}$ and the number of gossip iterations $L = \left\lceil \frac{2 \ln k}{\gamma \rho} \right\rceil$ for a constant \( k > 2 \). Denote by \(\Xb_t =[\Xb^{(1)}_{t},\ldots,\Xb^{(m)}_{t}]\) and let \(\bXb_t = \frac{1}{m}\Xb_t \mathbf{1}\) be the column-wise average.  If $\E[\sum_{i=1}^m\|\g_t^{(i)}\|_2^2]\le Q$ for all $t\in[T]$, with the inputs formed before the current consensus call, then the per-round communication error satisfies
	\begin{equation} \label{eqn:lem:choco-gossip}
		\E \left[ \| \Xb_t - \bXb_{t} \mathbf{1}^{\top} \|_F^2 + \|\Xb_t - \hXb_{t} \|_F^2 \right]
		\leq \frac{4Q}{k-2}, \qquad\forall t\in[T].
	\end{equation}
\end{lemma}

\subsection{The Proof of Theorem~\ref{thm:Semi_Dec_GDRO}}  
Recall the time-averaged local output $\wb_T^{(i)} = \frac{1}{T} \sum_{t=1}^T \w_t^{(i)}$, and define the time-averaged shared dual vector $\qb_T = \frac{1}{T} \sum_{t=1}^T \q_t$.
For any worker $i\in[m]$, we decompose the optimization error at $\wb_T^{(i)}$ as
\begin{equation}\label{eqn:proof:semi_dec:error_decom_1}
	\begin{split}
		&\varepsilon_{\phi}(\wb_T^{(i)})
		= \max_{j\in[m]}R_j(\wb_T^{(i)})-\min_{\w\in\mathcal{W}}\max_{j\in[m]}R_j(\w)=\max_{j\in[m]}R_j(\wb_T^{(i)})-\max_{j\in[m]}R_j(\w^*)\\
		= & \max_{\q\in \Delta_m}  \phi(\wb_T^{(i)},\q)- \max_{\q\in \Delta_m}  \phi(\w^*,\q) 
		\leq \max_{\q\in \Delta_m}  \phi(\wb_T^{(i)},\q)- \phi\left(\w^*,\qb_T\right)  \\
		\leq & \frac{1}{T} \max_{\q\in \Delta_{m}} \sum_{t=1}^T \phi(\w_t^{(i)},\q)- \frac{1}{T} \sum_{t=1}^T \phi(\w^*,\q_t) \\
		\leq & \frac{1}{T} \left(\max_{\q\in \Delta_{m}}  \sum_{t=1}^T \phi(\w_t^{(i)},\q)- \sum_{t=1}^T \phi(\w_t^{(i)},\q_t)\right) + \frac{1}{T} \left(\sum_{t=1}^T \phi(\w_t^{(i)},\q_t) - \sum_{t=1}^T\phi(\w^*,\q_t)\right).
	\end{split}
\end{equation}
As introduced in Section~\ref{sec:semi_dec_GDRO}, we rewrite 
$$\phi(\w_t^{(i)},\q)=\sum_{j=1}^m q_{j}R_j(\w_t^{(i)})=\sum_{j=1}^m q_{j}R_j(\w_t^{(j)})+\sum_{j=1}^m q_{j}(R_j(\w_t^{(i)})-R_j(\w_t^{(j)})).$$ We denote $j(i)=\argmax_{j\in[m]} \sum_{t=1}^T R_j(\w_t^{(i)})$ and obtain
\begin{equation} \label{eqn:proof:semi_dec_error_q1}
	\begin{split}
		&\max_{\q\in \Delta_{m}}  \sum_{t=1}^T \phi(\w_t^{(i)},\q)- \sum_{t=1}^T \phi(\w_t^{(i)},\q_t) \\
		=& \max_{\q\in \Delta_{m}}  \sum_{t=1}^T \sum_{j=1}^m q_j R_j(\w_t^{(i)}) - \sum_{t=1}^T\sum_{j=1}^m q_{t,j} R_j(\w_t^{(i)}) \\
        =& \sum_{t=1}^T R_{j(i)}(\w_t^{(i)}) - \sum_{t=1}^T\sum_{j=1}^m q_{t,j} R_j(\w_t^{(i)}) \\
		=&\sum_{t=1}^T \left(R_{j(i)}(\w_t^{(i)}) - R_{j(i)}(\w^{(j(i))}_t)\right) + \sum_{t=1}^T \left( R_{j(i)}(\w^{(j(i))}_t) - \sum_{j=1}^m q_{t,j} R_j(\w^{(j)}_t) \right) \\
		&+ \sum_{t=1}^T \sum_{j=1}^m q_{t,j} \left(R_j(\w^{(j)}_t) - R_j(\w_t^{(i)})\right) \\
		\overset{(\ref{eqn:ass:value_gradient})}{\leq} & G\sum_{t=1}^T\left(\|\w_t^{(i)}-\w^{(j(i))}_t\|_2+\sum_{j=1}^m q_{t,j}\|\w_t^{(i)}-\w^{(j)}_t\|_2\right) \\
        &+ \sum_{t=1}^T \left(R_{j(i)}(\w^{(j(i))}_t) - \sum_{j=1}^m q_{t,j} R_j(\w^{(j)}_t)\right).
	\end{split}
\end{equation}
where $s_{t,j}^{(i)} = 1 - R_j(\w_t^{(i)}) = 1 - \E_{\z \sim \P_j}[\ell(\w_t^{(i)};\z)] \overset{(\ref{eqn:ass:value_gradient})}{\in}  [0,1]$.
Substituting~(\ref{eqn:proof:semi_dec_error_q1}) into~(\ref{eqn:proof:semi_dec:error_decom_1}), we obtain
\begin{equation}  \label{eqn:proof:semi_dec:error_decom_2}
	\begin{split}
		\varepsilon_{\phi}(\wb_T^{(i)}) \leq & \frac{1}{T} \underbrace{\max_{k\in[m]}\sum_{t=1}^T\left(\sum_{j=1}^m q_{t,j} s_{t,j}^{(j)} - s_{t,k}^{(k)}\right)}_{\Regq} 
		+ \frac{1}{T} \underbrace{\sum_{t=1}^T\sum_{j=1}^m \left(q_{t,j}R_j(\w_t^{(i)}) - q_{t,j}R_j(\w^*)\right)}_{\Regw^{(i)}(\w^*)}\\
		&+\frac{1}{T}\underbrace{G\sum_{t=1}^T\left(\|\w_t^{(i)}-\w^{(j(i))}_t\|_2+\sum_{j=1}^m q_{t,j}\|\w_t^{(i)}-\w^{(j)}_t\|_2\right)}_{\mathtt{term(dis-semi)}}
	\end{split}
\end{equation}

\subsubsection{Regret $\Regw^{(i)}$} \label{sec:proof:semi_decGDRO:w}
For later reuse of the analysis, in SemiDec we write $q_{t,i}^{(i)}=q_{t,i}$, the weight sent by the server to worker $i$, for all $i\in[m]$ and $t\in[T]$.

Let $\bg_{w,t}= \frac{1}{m}\sum_{i=1}^m \g_{w,t}^{(i)}$. 
We introduce a virtual decision sequence $\{\tilde{\w}_t\}_{t\ge 1}$ initialized by $\tilde{\w}_{1}=\argmin_{\w \in \W} \nu_w(\w)$ and updated by
\begin{equation}  \label{eqn:proof:dec:virtual:w}
	\tilde{\w}_{t+1}=\argmin_{\w \in \W}\left\{ \left\langle\eta_w\sum_{\tau=1}^t\bg_{w,\tau},\w \right\rangle + \nu_w(\w) \right\} \quad\forall t\in[T].
\end{equation}
We characterize the deviation between each worker's local decision and the virtual decision as $\epswdist = \max_{i\in[m]}\bigl\|\w_t^{(i)}-\tilde{\w}_{t}\bigr\|_2$.

For a fixed worker $i\in[m]$, any distribution index $j\in[m]$ and any fixed point $\u\in\W$, we have 
\begin{equation}  \label{eqn:proof:semi_dec:reg_0}
	\begin{split}
		&R_j(\w_t^{(i)})-R_j(\u)= R_j(\w_t^{(j)})-R_j(\u) + R_j(\w_t^{(i)}) - R_j(\w_t^{(j)}) \\
		=&R_j(\w_t^{(j)})-R_j(\u) + \E_{\z \sim \P_j}[\ell(\w_t^{(i)};\z)]-\E_{\z \sim \P_j}[\ell(\w_t^{(j)};\z)] \\
		\overset{(\ref{eqn:ass:value_gradient})}{\leq} &R_j(\w_t^{(j)})-R_j(\u) + G\|\w_t^{(i)}-\w_t^{(j)}\|
		\leq \langle\nabla R_j(\w_t^{(j)}),\w_t^{(j)}-\u\rangle + G\|\w_t^{(i)}-\w_t^{(j)}\|\\
		\leq & \langle \nabla\ell(\w_t^{(j)};\z_t^{(j)}),\tilde{\w}_t-\u\rangle + \langle \nabla\ell(\w_t^{(j)};\z_t^{(j)}),\w_t^{(j)}-\tilde{\w}_t\rangle \\
		&+ \langle\nabla R_j(\w_t^{(j)})-\nabla\ell(\w_t^{(j)};\z_t^{(j)}),\w_t^{(j)}-\u\rangle + G\|\w_t^{(i)}-\w_t^{(j)}\|\\
		\leq & \langle \nabla\ell(\w_t^{(j)};\z_t^{(j)}),\tilde{\w}_t-\u\rangle + \langle\nabla R_j(\w_t^{(j)})-\nabla\ell(\w_t^{(j)};\z_t^{(j)}),\w_t^{(j)}-\u\rangle \\
		&+ G\left( \|\tilde{\w}_t-\w_t^{(j)}\|+\|\w_t^{(i)}-\w_t^{(j)}\|\right)\\
		\leq & \langle \nabla\ell(\w_t^{(j)};\z_t^{(j)}),\tilde{\w}_t-\u\rangle + \langle\nabla R_j(\w_t^{(j)})-\nabla\ell(\w_t^{(j)};\z_t^{(j)}),\w_t^{(j)}-\u\rangle
		+ 3G\epswdist.
	\end{split}
\end{equation}
Multiplying both sides of (\ref{eqn:proof:semi_dec:reg_0}) by $q_{t,j}^{(j)}$, using $q_{t,j}^{(j)}\nabla\ell(\w_t^{(j)};\z_t^{(j)})=\g_{w,t}^{(j)}$, summing over $j\in[m]$ and $t\in[T]$, and taking expectations, we obtain
\begin{equation}  \label{eqn:proof:semi_dec_reg_w:1}
	\begin{split}
		\E\left[\Regw^{(i)}(\u)\right] =& \E\left[\sum_{t=1}^T \sum_{j=1}^m q_{t,j}^{(j)} R_j(\w_t^{(i)}) - \sum_{t=1}^T\sum_{j=1}^m q_{t,j}^{(j)} R_j(\u)\right]   \\
		\leq & m\sum_{t=1}^T \E\left[\langle \bg_{w,t}, \tilde{\w}_{t}-\u\rangle\right] + 3G\sum_{t=1}^T\E[\epswdist].
	\end{split}
\end{equation}
For the first term on the right-hand side of~(\ref{eqn:proof:semi_dec_reg_w:1}), using Lemma~\ref{lem:FTRL_Orabona} with the update rule of $\tilde{\w}_{t+1}$ in (\ref{eqn:proof:dec:virtual:w}), we have
\begin{equation} \label{eqn:proof:semi_dec:w:1}
	\begin{split}
		\sum_{t=1}^T \langle  \bg_{w,t},\tilde{\w}_{t}-\u \rangle 
		\leq \ & \frac{\nu_w(\u) - \min_{\w\in \W}\nu_w(\w) }{\eta_w} + \frac{\eta_w}{2}\sum_{t=1}^T \|\bg_{w,t}\|_2^2 \\
		\leq \ & \frac{D^2}{\eta_w} + \frac{\eta_w}{2}\sum_{t=1}^T \left\|\frac{1}{m}\sum_{i=1}^m q_{t,i}^{(i)}\cdot\nabla\ell(\w_t^{(i)};\z_t^{(i)})\right\|_2^2\\
		\overset{(\ref{eqn:ass:value_gradient})}{\leq} & \frac{D^2}{\eta_w} + \frac{\eta_w TG^2}{2m^2},
	\end{split}
\end{equation}
where the second inequality follows from Assumption~\ref{ass:domain_bound}, and the last inequality uses the fact that, in SemiDec, the centralized update of $\q_t$ guarantees that $\sum_{i=1}^m q_{t,i}^{(i)} = \sum_{i=1}^m q_{t,i} = 1$.

Combining~(\ref{eqn:proof:semi_dec_reg_w:1}) and~(\ref{eqn:proof:semi_dec:w:1}) and setting $\eta_w=\frac{\sqrt{2}mD}{\sqrt{T}G}$, we obtain
\begin{equation} \label{eqn:proof:semi_dec:w:res}
	\begin{split}
		\E\left[\Regw^{(i)}(\u)\right] =& \E\left[\sum_{t=1}^T \sum_{j=1}^m q_{t,j}^{(j)} R_j(\w_t^{(i)}) - \sum_{t=1}^T\sum_{j=1}^m q_{t,j}^{(j)} R_j(\u)\right]   \\
		\leq & \frac{mD^2}{\eta_w} + \frac{\eta_w TG^2}{2m} + 3G\sum_{t=1}^T\E[\epswdist]
        =  \sqrt{2T}DG+ 3G\sum_{t=1}^T \E[\epswdist],
	\end{split}
\end{equation}

\subsubsection{Regret $\Regq$}   \label{sec:proof:semi_decGDRO:q}
In this part, we derive a bound on the regret $\Regq$ in the semi-decentralized setting. Recall that the server solves a \emph{prediction with expert advice} (PEA) problem and adopts the Hedge algorithm~\citep{Hedge_1997}. Following~\citet[Lemma 1]{Bai2025FSQ}, let $\hats_{t,j}^{(i)} = 1 - \ell(\w_t^{(i)};\z^{(i)}_t), \ \tilde{\s}_t=[\hats_{t,1}^{(1)};\ldots;\hats_{t,m}^{(m)}]$, and $\s_t=(s_{t,1}^{(1)},\ldots,s_{t,m}^{(m)})^\top$. We decompose $\Regq$ as
\begin{equation} \label{eqn:decom:error:q}
	\begin{split}
		\Regq 
		\leq \underbrace{\sum_{t=1}^T \inner{\q_t}{\s_t-\tilde{\s}_t}}_{\mathtt{term (q_1)}} + \underbrace{\sum_{t=1}^T \inner{\q_t}{\tildesb_t} - \min_{i\in[m]} \sum_{t=1}^T \hats_{t,i}^{(i)}}_{\mathtt{term (q_2)}} + \underbrace{\sum_{t=1}^T \hats_{t,k^*}^{(k^*)} - \sum_{t=1}^T s_{t,k^*}^{(k^*)}}_{\mathtt{term (q_3)}}.
	\end{split}
\end{equation}
where $k^* \in \argmin_{k\in [m]} \sum_{t=1}^T s_{t,k}^{(k)}$.

Let $\F_{t-1}$ denote the history before round $t$. For $\mathtt{term (q_1)}$, since $\q_t$ is $\F_{t-1}$-measurable and
$\E[\tilde{\s}_t\mid \F_{t-1}]=\s_t$, we have $\E[\mathtt{term (q_1)}]=0$.

For $\mathtt{term (q_2)}$, by Lemma~\ref{lem:FTRL_Orabona:local_norm} and setting
$\eta_q=\sqrt{\frac{2\ln m}{T}}$, we can derive for each $k\in[m]$,
\begin{equation}  \label{eqn:proof:q1}
\begin{split}
\sum_{t=1}^T \left(\inner{\q_t}{\tildesb_t}-\hats_{t,k}^{(k)}\right)
&\leq \frac{\ln m}{\eta_{q}}
+ \frac{\eta_{q}}{2}\sum_{t=1}^T\sum_{k=1}^{m}q_{t,k}(\hats_{t,k}^{(k)})^{2}\\
&\leq \frac{\ln m}{\eta_{q}}+ \frac{T\eta_{q}}{2}
= \sqrt{2T\ln m}.
\end{split}
\end{equation}
Thus, we obtain $\mathtt{term (q_2)}\le \sqrt{2T\ln m}$.

For $\mathtt{term (q_3)}$, define $X_{t,k}=\hats_{t,k}^{(k)}-s_{t,k}^{(k)}$. Then $\E[X_{t,k}\mid \F_{t-1}]=0$. Moreover, by Assumption~\ref{ass:basic_pointwise},
$\hats_{t,k}^{(k)}\in[0,1]$, we obtain $|X_{t,k}|\leq1$. 
By Lemma~\ref{lem:azuma} and a union bound over \(i\in[m]\), with probability at least \(1-\delta\),
\[
\mathtt{term}(\q_3)\leq \max_{k\in[m]}\sum_{t=1}^T X_{t,k}\leq \sqrt{2T\ln\frac{m}{\delta}}.
\]
We now convert this high-probability bound into an expectation bound. Let \(Z=\frac{(\mathtt{term}(\q_3))_+^2}{2T}-\ln m\). The above inequality implies that, for any \(\delta\in(0,1)\),
\[
\Pr\left[Z>\ln\frac{1}{\delta}\right]=\Pr\left[\mathtt{term}(\q_3)>\sqrt{2T\ln\frac{m}{\delta}}\right]\leq \delta.
\]
Applying Lemma~\ref{lemma:high_pro_bound_to_Exp} to \(Z\), we get
\[
\E[Z]\leq \int_0^1 \frac{1}{\delta}\Pr\left[Z>\ln\frac{1}{\delta}\right]\mathrm d\delta\leq 1.
\]
Therefore, 
\begin{equation} \label{eqn:proof:q3}
    \E[\mathtt{term}(\q_3)]\leq \E[(\mathtt{term}(\q_3))_+]\leq \sqrt{\E[(\mathtt{term}(\q_3))_+^2]}\leq \sqrt{2T(\ln (m)+1)}.
\end{equation}

Taking expectations on (\ref{eqn:decom:error:q}), combining it with the three bounds and using $\ln m\ge\ln2>1/2$ and $\sqrt2+\sqrt6<4$, we obtain
\begin{equation} \label{eqn:decom:error:q:res}
\E\left[\Regq\right]\leq 4\sqrt{T\ln m}.
\end{equation}

\subsubsection{Discrepancy error term and results}
To control the discrepancy term, recalling that $\epswdist=\max_{i\in[m]}\|\w_t^{(i)}-\tilde{\w}_t\|_2$, we obtain
\begin{equation}\label{eqn:proof:estimator_error1}
	\begin{split}
		\mathtt{term(dis\text{-}semi)}
		=&\ G\sum_{t=1}^T\left(\|\w_t^{(i)}-\w^{(j(i))}_t\|_2+\sum_{j=1}^m q_{t,j}\|\w_t^{(i)}-\w^{(j)}_t\|_2\right) \\
		\leq&\ G\sum_{t=1}^T\left(2\epswdist+\sum_{j=1}^m q_{t,j}\cdot 2\epswdist\right)
		= 4G\sum_{t=1}^T \epswdist .
	\end{split}
\end{equation}

Substituting~(\ref{eqn:proof:semi_dec:w:res}), (\ref{eqn:decom:error:q:res}), (\ref{eqn:proof:estimator_error1}) into~(\ref{eqn:proof:semi_dec:error_decom_2}), setting $\eta_w=\frac{\sqrt{2}mD}{\sqrt{T}G}$, $\eta_q=\sqrt{\frac{2\ln m}{T}}$, and taking expectations yields
\begin{equation*}\label{eqn:proof:semi_dec:error_decom_3_pre}
	\begin{split}
		\E\left[\varepsilon_{\phi}(\wb_T^{(i)})\right]
		\leq\ & \frac{1}{T}\E[\Regq] + \frac{1}{T}\E\left[\Regw^{(i)}(\w^*)\right]
		+ \frac{1}{T}\E\left[\mathtt{term(dis\text{-}semi)}\right] \\
		\leq\ & 4\sqrt{\frac{\ln m}{T}}
		+ \sqrt{\frac{2}{T}}DG
		+ \frac{3G}{T}\sum_{t=1}^T \E[\epswdist]
		+ \frac{4G}{T}\sum_{t=1}^T \E[\epswdist] \\
		=\ & 4\sqrt{\frac{\ln m}{T}}
		+ \sqrt{\frac{2}{T}}DG
		+ \frac{7G}{T}\sum_{t=1}^T \E[\epswdist].
	\end{split}
\end{equation*}
For any $i\in[m]$ and $t>1$, by Lemma~\ref{lem:Projection} and the update rule~(\ref{eqn:update:wt:dec}) and (\ref{eqn:proof:dec:virtual:w}), we have
\begin{equation}\label{eqn:proof:semi_dec:epswdist_decom}
		\|\w_t^{(i)}-\tilde{\w}_t\|_2
		\overset{(\ref{eqn:lemma:projection})}{\leq}\ 
		\eta_w \left\|\sum_{\tau=1}^{t-1}\bg_{w,\tau}-\Fb_{t-1}^{(i)}\right\|_2 
		\leq\ 
		\eta_w \left(
		\left\|\Fb_{t-1}^{(i)}-\bFb_{t-1}\right\|_2
		+\left\|\bFb_{t-1}-\sum_{\tau=1}^{t-1}\bg_{w,\tau}\right\|_2
		\right),
\end{equation}
If the consensus routine guarantees $\bFb_{t-1}=\sum_{\tau=1}^{t-1}\bg_{w,\tau}$ and $\E[\max_{r\in[m]}\|\Fb_{t-1}^{(r)}-\bFb_{t-1}\|_2]\leq\varepsilon_c$ for all $t\in[T]$, then for any worker $i\in[m]$ we obtain
\begin{equation}\label{eqn:proof:semi_dec:thm_res1}
\begin{aligned}
\E\left[\varepsilon_{\phi}(\wb_T^{(i)})\right]
\le&
4\sqrt{\frac{\ln m}{T}}
+\sqrt{\frac{2}{T}}DG
+\frac{7\sqrt{2}mD}{T\sqrt{T}}
\sum_{t=1}^T
\E\left[\max_{r\in[m]}\left\|\Fb_{t-1}^{(r)}-\bFb_{t-1}\right\|_2\right]\\
\leq & 4\sqrt{\frac{\ln m}{T}}
+\sqrt{\frac{2}{T}}DG
+\frac{7\sqrt{2}mD}{\sqrt{T}}\varepsilon_c
\end{aligned}
\end{equation}

\subsection{The Proof of Corollary~\ref{cor:Semi_Dec_AG}}   \label{app:proof:Semi_Dec_AG}
When adopting the accelerated gossip update~(\ref{lem:acc_gossip_update}), we set $\Fb_t^{(i)} = \Xb_t^{(i)}$ and $\bFb_t = \bXb_t = \frac{1}{m}\Xb_t\mathbf{1}$.
By Lemma~\ref{lem:iter_bound:ag_gossip}, the average is preserved, i.e., $\bFb_t = \bFb_{t-1} + \bg_t$, where $\bg_t = \frac{1}{m}\sum_{i=1}^m \g_t^{(i)}$.
Since $\bFb_0=\mathbf{0}$, we obtain by induction that $\bFb_{t-1} = \sum_{\tau=1}^{t-1}\bg_{w,\tau}$, and hence the second term in~(\ref{eqn:proof:semi_dec:epswdist_decom}) vanishes. Let $a_m=\sqrt{(m-1)/m}$. For centered columns $\sum_i v_i=0$, Cauchy--Schwarz gives $\|v_i\|_2^2\le(m-1)\sum_{j\ne i}\|v_j\|_2^2$, and thus $\max_i\|v_i\|_2\le a_m(\sum_i\|v_i\|_2^2)^{1/2}$. Applying this inequality and the per-round communication error bound~(\ref{eq:lem:ag_gossip}) in Lemma~\ref{lem:iter_bound:ag_gossip}, we obtain
\begin{equation}\label{eqn:proof:semi_dec:epswdist_AG_final}
\|\Fb_{t-1}^{(i)}-\bFb_{t-1}\|_2\leq \frac{a_mG}{\sqrt{k}-1},\qquad \forall t\in[T],
\end{equation}
since $\sum_{i=1}^m \|\g_t^{(i)}\|_2^2=\sum_{i=1}^m \|q_{t,i}\cdot\nabla\ell(\w_t^{(i)};\z_t^{(i)})\|_2^2 \le G^2$ for all $t\in[T]$. 
Combining (\ref{eqn:proof:semi_dec:epswdist_AG_final}) with (\ref{eqn:proof:semi_dec:thm_res1}), and setting $k=(\frac{m}{\sqrt{\ln m}}+1)^2$, we obtain
\begin{equation}\label{eqn:proof:semi_dec:final_bound}
\begin{split}
    \E\left[\varepsilon_{\phi}(\wb_T^{(i)})\right]
\leq4\sqrt{\frac{\ln m}{T}}+ \sqrt{\frac{2}{T}}DG+ 7\sqrt{2}a_mDG\sqrt{\frac{\ln m}{T}} \\
\leq (8\sqrt{2}DG+4)\sqrt{\frac{\ln m}{T}}
= O\left(\sqrt{\frac{\log m}{T}}\right).
\end{split}
\end{equation}
The last inequality uses $(\ln m)^{-1/2}+7a_m\le8$: for $m\ge3$, use $\ln m>1$ and $a_m\le1$; for $m=2$, use $\ln2>1/2$ and $a_2=1/\sqrt2$.

\subsection{The Proof of Corollary~\ref{cor:Semi_Dec_CG}} \label{app:proof:Semi_Dec_CG}
We first show that the second term in~(\ref{eqn:proof:semi_dec:epswdist_decom}) vanishes, i.e., $\bFb_{t}-\sum_{\tau=1}^{t}\bg_{w,\tau}
=\frac{1}{m}\sum_{i=1}^{m}\Fb_{t}^{(i)}-\frac{1}{m}\sum_{\tau=1}^{t}\sum_{i=1}^{m}\g_{w,\tau}^{(i)}
=\mathbf{0}_d$.
This follows from the average-preserving property of Algorithm~\ref{alg:choco-gossip}. In particular, for the inner gossip iterates, we have
\begin{equation}\label{eqn:proof:dec:gossip:the_same}
	\begin{split}
		\frac{1}{m}\sum_{i=1}^m \x_{L+1}^{(i)}
		&= \frac{1}{m}\sum_{i=1}^m\left(\x_{L}^{(i)} + \gamma\sum_{j\in \mathcal{N}_i}P_{i,j} \big(\hat{\x}_{L+1}^{(j)}-\hat{\x}_{L+1}^{(i)}\big)\right) \\
		&= \frac{1}{m}\sum_{i=1}^m \x_{L}^{(i)}
		= \cdots
		= \frac{1}{m}\sum_{i=1}^m \x_{1}^{(i)},
	\end{split}
\end{equation}
where the second equality uses Assumption~\ref{ass:graph}.
Applying this property to the outer update of $\Fb_t^{(i)}$ yields
\begin{equation}\label{eqn:proof:dec:w:the_same}
		\bFb_{t}
		=\frac{1}{m}\sum_{i=1}^m \Fb_t^{(i)}
		= \frac{1}{m}\sum_{i=1}^m\left(\Fb^{(i)}_{t-1} + \g_{w,t}^{(i)}\right) 
		= \bFb_{t-1} + \bg_{w,t}
		= \sum_{\tau=1}^{t}\bg_{w,\tau},
\end{equation}
where the last equality uses $\bFb_0=\mathbf{0}_d$.
Consequently, we obtain $\bFb_{t-1}=\sum_{\tau=1}^{t-1}\bg_{w,\tau}$. 
Applying the per-round communication error bound~(\ref{eqn:lem:choco-gossip}) in Lemma~\ref{lem:semi_dec:choco-gossip} and Jensen's inequality yields for all $t\in[T]$
\begin{equation}\label{eqn:proof:semi_dec:epswdist_CG_final}
\E\left[\|\Fb_{t-1}-\bFb_{t-1}\mathbf{1}^\top\|_F\right]
\le  \sqrt{\E\left[\|\Fb_{t-1}-\bFb_{t-1}\mathbf{1}^\top\|_F^2\right]}
\le  \sqrt{\frac{4Q}{k-2}}
= 2\frac{G}{\sqrt{k-2}}.
\end{equation}
since $\sum_{i=1}^m \|\g_t^{(i)}\|_2^2
= \sum_{i=1}^m \|q_{t,i}\nabla\ell(\w_t^{(i)};\z_t^{(i)})\|_2^2
\le G^2$ for all $t\in[T]$.

Using the centered-column inequality from Appendix~\ref{app:proof:Semi_Dec_AG}, with $a_m=\sqrt{(m-1)/m}$, we may take $\varepsilon_c=2a_mG/\sqrt{k-2}$. Combining~(\ref{eqn:proof:semi_dec:epswdist_CG_final}) with~(\ref{eqn:proof:semi_dec:thm_res1}) and setting $k=\frac{m^2}{\ln m}+2$, we obtain
\begin{equation}\label{eqn:proof:semi_dec:final_bound_CG}
\begin{split}
&\E\left[\varepsilon_{\phi}(\wb_T^{(i)})\right]\\
\le&
4\sqrt{\frac{\ln m}{T}}
+ \sqrt{\frac{2}{T}}DG
+ 14\sqrt{2}a_mDG\sqrt{\frac{\ln m}{T}}
\leq (15\sqrt{2}DG+4)\sqrt{\frac{\ln m}{T}}
= O\left(\sqrt{\frac{\log m}{T}}\right).
\end{split}
\end{equation}
The last inequality uses $(\ln m)^{-1/2}+14a_m\le15$, which follows from $\ln m>1$ for $m\ge3$, and from $\ln2>1/2$ and $a_2=1/\sqrt2$ for $m=2$.

\subsection{The Proof of Theorem~\ref{thm:Dec_GDRO_modular}} \label{app:proof:Dec_GDRO_modular}
Let $\qb_T^{(i)}=T^{-1}\sum_{t=1}^T\q_t^{(i)}$. We decompose the optimization error at $\wb_T^{(i)}$ as
\begin{equation}\label{eqn:proof:dec:error_decom_1}
	\begin{split}
		&\varepsilon_{\phi}(\wb_T^{(i)})
		= \max_{j\in[m]}R_j(\wb_T^{(i)})-\min_{\w\in\mathcal{W}}\max_{j\in[m]}R_j(\w)=\max_{j\in[m]}R_j(\wb_T^{(i)})-\max_{j\in[m]}R_j(\w^*)\\
		= & \max_{\q\in \Delta_m}  \phi(\wb_T^{(i)},\q)- \max_{\q\in \Delta_m}  \phi(\w^*,\q) 
		\leq \max_{\q\in \Delta_m}  \phi(\wb_T^{(i)},\q)- \phi\left(\w^*,\qb_T^{(i)}\right)  \\
		\leq & \frac{1}{T} \max_{\q\in \Delta_{m}} \sum_{t=1}^T \phi(\w_t^{(i)},\q)- \frac{1}{T} \sum_{t=1}^T \phi(\w^*,\q_t^{(i)}) \\
		\leq & \frac{1}{T} \left(\max_{\q\in \Delta_{m}}  \sum_{t=1}^T \phi(\w_t^{(i)},\q)- \sum_{t=1}^T \phi(\w_t^{(i)},\q_t^{(i)})\right) + \frac{1}{T} \left(\sum_{t=1}^T \phi(\w_t^{(i)},\q_t^{(i)}) - \sum_{t=1}^T\phi(\w^*,\q_t^{(i)})\right) 
	\end{split}
\end{equation}
In the decentralized setting, there is no central server to compute a shared dual sequence $\{\q_t\}$. 
Consequently, the above decomposition is expressed in terms of regret defined with respect to each worker's local dual iterates $\{\q_t^{(i)}\}$.

Recall that $j(i)=\argmax_{j\in[m]} \sum_{t=1}^T R_j(\w_t^{(i)})$ and $s_{t,j}^{(i)} = 1 - R_j(\w_t^{(i)})$. For the first term, following~(\ref{eqn:proof:semi_dec_error_q1}), we have for all $i\in[m]$,
\begin{equation} \label{eqn:proof:dec_error_q1}
	\begin{split}
		&\max_{\q\in \Delta_{m}}  \sum_{t=1}^T \phi(\w_t^{(i)},\q)- \sum_{t=1}^T \phi(\w_t^{(i)},\q_t^{(i)}) \\
		\leq &G\sum_{t=1}^T\left(\|\w_t^{(i)}-\w^{(j(i))}_t\|_2+\sum_{j=1}^m q_{t,j}^{(i)}\|\w_t^{(i)}-\w^{(j)}_t\|_2\right) + \sum_{t=1}^T\sum_{j=1}^m q_{t,j}^{(i)} s_{t,j}^{(j)} - \sum_{t=1}^T s_{t,j(i)}^{(j(i))}
	\end{split}
\end{equation}
For the second term, we define $\hat{\q}_{t}=[q^{(1)}_{t,1},\ldots,q^{(m)}_{t,m}]^{\top}$ and have 
\begin{equation}\label{eqn:proof:dec_error_w1}
	\begin{split}
		&\sum_{t=1}^T \phi(\w_t^{(i)},\q_t^{(i)}) - \sum_{t=1}^T\phi(\w^*,\q_t^{(i)})\\
		= &\sum_{t=1}^T \sum_{j=1}^m \left(q^{(i)}_{t,j}R_j(\w_t^{(i)}) -q^{(j)}_{t,j}R_j(\w_t^{(i)})\right) + \sum_{t=1}^T \sum_{j=1}^m \left(q^{(j)}_{t,j}R_j(\w_t^{(i)}) -q^{(j)}_{t,j}R_j(\w^*)\right) \\
		&+ \sum_{t=1}^T \sum_{j=1}^m \left(q^{(j)}_{t,j}R_j(\w^*) -q^{(i)}_{t,j}R_j(\w^*)\right) \\
		\overset{(\ref{eqn:ass:value_gradient})}{\leq} &2\sum_{t=1}^T \|\hat{\q}_{t}-\q_t^{(i)}\|_1 + \sum_{t=1}^T \sum_{j=1}^m \left(q^{(j)}_{t,j}R_j(\w_t^{(i)}) -q^{(j)}_{t,j}R_j(\w^*)\right)
	\end{split}
\end{equation}
Combining~(\ref{eqn:proof:dec_error_q1}) and~(\ref{eqn:proof:dec_error_w1}) with the decomposition in~(\ref{eqn:proof:dec:error_decom_1}), we obtain the following bound:
\begin{equation}\label{eqn:proof:dec:error_decom_final_dec}
	\begin{split}
		\varepsilon_{\phi}(\wb_T^{(i)})  
		\leq & \frac{1}{T}\underbrace{\sum_{t=1}^T \sum_{j=1}^m \left(q^{(j)}_{t,j}R_j(\w_t^{(i)}) -q^{(j)}_{t,j}R_j(\w^*)\right)}_{\Regw^{(i)}(\w^*)}+ \frac{1}{T}\underbrace{\sum_{t=1}^T\left(\sum_{j=1}^m q_{t,j}^{(i)} s_{t,j}^{(j)} - s_{t,j(i)}^{(j(i))}\right)}_{\Regq^{(i)}}\\
		& +\frac{1}{T} \underbrace{\sum_{t=1}^T\left(2\|\hat{\q}_{t}-\q_t^{(i)}\|_1+G\|\w_t^{(i)}-\w^{(j(i))}_t\|_2+G\sum_{j=1}^m q_{t,j}^{(i)}\|\w_t^{(i)}-\w^{(j)}_t\|_2\right)}_{{\mathtt{term(dis-dec)}}}
	\end{split}
\end{equation}
Let $\bg_{q,t}=m^{-1}\sum_{i=1}^m\g_{q,t}^{(i)}$. We define the virtual iterate $\tilde{\w}_{t+1}$ as in~(\ref{eqn:proof:dec:virtual:w}) and 
\begin{equation}  \label{eqn:proof:dec:virtual:q}
	\tilde{\q}_{t+1}=\argmin_{\q \in \Delta_m}\left\{ \left\langle\eta_q\sum_{\tau=1}^t\bg_{q,\tau},\q \right\rangle + \nu_q\left(\q\right) \right\},
\end{equation}
with the initialization $\tilde{\q}_1=\argmin_{\q \in \Delta_m}\nu_q\left(\q\right)$.
To streamline the subsequent bounds, we introduce the following discrepancy measures
\[
\epswdist = \max_{i\in[m]}\bigl\|\w_t^{(i)}-\tilde{\w}_t\bigr\|_2,
\qquad
\epsqdist = \max_{i\in[m]}\bigl\|\q_t^{(i)}-\tilde{\q}_t\bigr\|_2,
\qquad t\in[T].
\]
These quantities capture, respectively, the disagreement of the primal and dual iterates
around the corresponding virtual sequences. In the following, we analyze each of these components in~(\ref{eqn:proof:dec:error_decom_final_dec}).

\subsubsection{Regret $\Regw^{(i)}$}   \label{sec:proof:DecGDRO:w}
The analysis of the local regret $\Regw^{(i)}$ in the decentralized setting closely follows the SemiDec analysis in Section~\ref{sec:proof:semi_decGDRO:w}. The main difference is that the regret uses the locally maintained dual weights $q_{t,i}^{(i)}$ for all $i\in[m]$, rather than the corresponding weights from the shared dual vector $\q_t$.

Note that, in the decentralized setting, the averaged gradient is written as
$\bg_{w,t}= \frac{1}{m}\sum_{i=1}^m \g_{w,t}^{(i)}=\frac{1}{m}\sum_{i=1}^m q_{t,i}^{(i)}\cdot\nabla\ell(\w_t^{(i)};\z_t^{(i)})$. Following~(\ref{eqn:proof:semi_dec_reg_w:1}), we obtain
\begin{equation}  \label{eqn:proof:dec_reg_w:1}
\begin{split}
\E[\Regw^{(i)}(\u)]
&\le m\sum_{t=1}^T\E[\langle\bg_{w,t},\tilde{\w}_t-\u\rangle]\\
&\quad+G\sum_{t=1}^T\E\left[\sum_{j=1}^m q_{t,j}^{(j)}
\left(\|\tilde{\w}_t-\w_t^{(j)}\|_2+\|\w_t^{(i)}-\w_t^{(j)}\|_2\right)\right].
\end{split}
\end{equation}
For the first term on the right-hand side of~(\ref{eqn:proof:dec_reg_w:1}), we proceed as in the proof of~(\ref{eqn:proof:semi_dec:w:1}). 
In the decentralized setting, the averaged gradient can be rewritten as
\[
\bg_{w,t}=\frac{1}{m}\sum_{i=1}^m q_{t,i}^{(i)}\nabla\ell(\w_t^{(i)};\z_t^{(i)})
=\frac{1}{m}\sum_{i=1}^m \bigl(q_{t,i}^{(i)}-\tilde{q}_{t,i}\bigr)\nabla\ell(\w_t^{(i)};\z_t^{(i)})
+\frac{1}{m}\sum_{i=1}^m \tilde{q}_{t,i}\nabla\ell(\w_t^{(i)};\z_t^{(i)}).
\]
Under Assumption~\ref{ass:basic_pointwise}, we have
\begin{equation}\label{eqn:proof:dec:avg_grad_decomp}
\begin{split}
\|\bg_{w,t}\|_2
&\le \frac{1}{m}\sum_{i=1}^m |q_{t,i}^{(i)}-\tilde{q}_{t,i}|\cdot \|\nabla\ell(\w_t^{(i)};\z_t^{(i)})\|_2
+\left\|\frac{1}{m}\sum_{i=1}^m \tilde{q}_{t,i}\nabla\ell(\w_t^{(i)};\z_t^{(i)})\right\|_2 \\
&\le \frac{G}{m}\|\hat{\q}_t-\tilde{\q}_t\|_1+\frac{G}{m}\sum_{i=1}^m \tilde{q}_{t,i}
= \left(1+\|\hat{\q}_t-\tilde{\q}_t\|_1\right)\frac{G}{m}\leq  \left(1+m \epsqdist\right)\frac{G}{m},
\end{split}
\end{equation}
where the last inequality follows 
\begin{equation}  \label{eqn:eps_q_1to2}
\|\hat{\q}_{t}-\tilde{\q}_t\|_1= \sum_{i=1}^m |q_{t,i}^{(i)}-\tilde{q}_{t,i}|
\leq \sum_{i=1}^m \|\q_{t}^{(i)}-\tilde{\q}_t\|_2
\leq m \epsqdist.
\end{equation}
Consequently,
\begin{equation}\label{eqn:proof:dec:regw_term1_bound}
\begin{aligned}
m\sum_{t=1}^T \E\left[\langle \bg_{w,t},\tilde{\w}_{t}-\u\rangle\right]
&\le m\left(\frac{D^2}{\eta_w}+\frac{\eta_w}{2}\sum_{t=1}^T \E\left[\|\bg_{w,t}\|_2^2\right]\right)\\
&\le \frac{mD^2}{\eta_w}+\frac{\eta_w G^2}{2m}\sum_{t=1}^T \E\left[\left(1+m\epsqdist\right)^2\right].
\end{aligned}
\end{equation}

We now bound the last term on the right-hand side of~(\ref{eqn:proof:dec_reg_w:1}). By the definition of $\epswdist$, we have
\begin{equation}  \label{eqn:proof:dec:regw_term2_bound}
\begin{split}
    &\sum_{t=1}^T\sum_{j=1}^m q_{t,j}^{(j)}\left( \|\tilde{\w}_t-\w_t^{(j)}\|+\|\w_t^{(i)}-\w_t^{(j)}\|\right)\leq 3\sum_{t=1}^T\sum_{j=1}^m q_{t,j}^{(j)}\epswdist \\
   = &3\sum_{t=1}^T\sum_{j=1}^m \tilde{q}_{t,j}\epswdist + 3\sum_{t=1}^T\sum_{j=1}^m \left(q_{t,j}^{(j)}-\tilde{q}_{t,j}\right)\epswdist \\
   \leq& 3\sum_{t=1}^T\left(1+\|\hat{\q}_t-\tilde{\q}_t\|_1\right)\epswdist
   \leq 3\sum_{t=1}^T\left(1+m \epsqdist\right)\epswdist,
   \end{split}
\end{equation}
where the last inequality follows~(\ref{eqn:eps_q_1to2}).
Substituting~(\ref{eqn:proof:dec:regw_term1_bound}) and~(\ref{eqn:proof:dec:regw_term2_bound}) into~(\ref{eqn:proof:dec_reg_w:1}) and setting $\eta_w=\frac{\sqrt{2}mD}{\sqrt{T}G}$ yields
\begin{equation}\label{eqn:proof:dec:regw_final_bound}
\ \E\left[\Regw^{(i)}(\u)\right]
\le
\sqrt{\frac{T}{2}}DG
+\frac{DG}{\sqrt{2T}}\sum_{t=1}^T \E\left[\left(1+m\epsqdist\right)^2\right]
+3G\sum_{t=1}^T \E\left[\left(1+m\epsqdist\right)\epswdist\right].
\end{equation}

\subsubsection{Regret $\Regq^{(i)}$}   \label{sec:proof:DecGDRO:q}
Recall $s_{t,j}^{(i)} = 1 - R_j(\w_t^{(i)}), \ \hats_{t,j}^{(i)} = 1 - \ell(\w_t^{(i)};\z^{(i)}_t), \ \tilde{\s}_t=[\hats_{t,1}^{(1)};\ldots;\hats_{t,m}^{(m)}]$. We then decompose the regret \( \Regq^{(i)} \) as
\begin{equation} \label{eqn:decom:dec_error:q:exp}
	\begin{split}
		&\E\left[\Regq^{(i)}\right] 
        = \E\left[\sum_{t=1}^T\left(\sum_{j=1}^m q_{t,j}^{(i)} s_{t,j}^{(j)} - s_{t,j(i)}^{(j(i))}\right)\right]=\E\left[\sum_{t=1}^T \inner{\q_t^{(i)}}{\tilde{\s}_t} - \sum_{t=1}^T s_{t,j(i)}^{(j(i))}\right]\\
		=& \underbrace{\E\left[\sum_{t=1}^T \left(\inner{\q_t^{(i)}}{\tilde{\s}_t} - \inner{\tilde{\q}_{t}}{\tilde{\s}_t}\right)\right]}_{\mathtt{term (q_1)}}+\underbrace{\E\left[\sum_{t=1}^T\left(\inner{\tilde{\q}_{t}}{\tilde{\s}_t} - \hats_{t,j(i)}^{(j(i))}\right)\right]}_{\mathtt{term (q_2)}}\\
		&+\underbrace{\E\left[\sum_{t=1}^T\left(\hats_{t,j(i)}^{(j(i))}-s_{t,j(i)}^{(j(i))}\right)\right]}_{\mathtt{term (q_3)}}
	\end{split}
\end{equation}
For $\mathtt{term (q_1)}$, we have 
\begin{equation} \label{eqn:proof:dec:q:term1}
	\begin{split}
		\sum_{t=1}^T \left(\inner{\q_t^{(i)}}{\tilde{\s}_t} - \inner{\tilde{\q}_{t}}{\tilde{\s}_t}\right)
		\leq \sum_{t=1}^T \|\q_t^{(i)}-\tilde{\q}_{t}\|_2 \|\tilde{\s}_t\|_2 
		\overset{(\ref{eqn:ass:value_gradient})}{\leq}\sqrt{m} \sum_{t=1}^T \|\q_t^{(i)}-\tilde{\q}_{t}\|_2 \leq \sqrt{m} \sum_{t=1}^T \epsqdist.
	\end{split}
\end{equation}
Then, we analyze $\mathtt{term (q_2)}$. By Lemma~\ref{lem:FTRL_Orabona:local_norm} and the update rule~(\ref{eqn:proof:dec:virtual:q}), we obtain for all $i\in[m]$,
\begin{equation}  \label{eqn:proof:dec:q1}
	\begin{split}
		\sum_{t=1}^T \left(\inner{\tilde{\q}_{t}}{\bg_{q,t}}-\frac{1}{m}\hats_{t,i}^{(i)}\right)
		\leq \frac{\ln m}{\eta_{q}} + \frac{\eta_{q}}{2}\sum_{t=1}^T\sum_{i=1}^{m}\tilde{q}_{t,i}\left(\frac{1}{m}\hats_{t,i}^{(i)}\right)^{2}\overset{(\ref{eqn:ass:value_gradient})}{\leq} \frac{\ln m}{\eta_{q}}+ \frac{\eta_{q}T}{2m^2},
	\end{split}
\end{equation}
Noting that $\tilde{\s}_t=m\bg_{q,t}$, we have 
\begin{equation} \label{eqn:proof:dec:q:term2_res}
	\E\left[\sum_{t=1}^T\left(\inner{\tilde{\q}_{t}}{\tilde{\s}_t} - \hats_{t,j(i)}^{(j(i))}\right)\right]\overset{(\ref{eqn:proof:dec:q1})}{\leq} \frac{m\ln m}{\eta_{q}}+ \frac{\eta_{q}T}{2m}.
\end{equation}
For $\mathtt{term (q_3)}$, we apply the argument for~(\ref{eqn:proof:q3}) to the data-dependent index $j(i)$ and obtain
\begin{equation} \label{eqn:proof:dec:q3_res}
    \mathtt{term}(q_3)\leq \sqrt{2T(\ln (m)+1)}\leq \sqrt6\,\sqrt{T\ln m}.
\end{equation}
Here we use $\ln m\ge\ln2>1/2$. Combining~(\ref{eqn:decom:dec_error:q:exp}), (\ref{eqn:proof:dec:q:term1}), (\ref{eqn:proof:dec:q:term2_res})--(\ref{eqn:proof:dec:q3_res}), using $\sqrt2+\sqrt6<4$, and setting $\eta_q=m\sqrt{\frac{2\ln m}{T}}$, we obtain
\begin{equation}\label{eqn:proof:dec:q:final_bound}
\E\left[\Regq^{(i)}\right] \le \sqrt{m}\sum_{t=1}^T \E[\epsqdist]+4\sqrt{T\ln m}.
\end{equation}

\subsubsection{Discrepancy error term and results}
We now bound the error caused by the disagreement among local models and between each local dual vector and the proxy vector, captured by
\begin{equation}
    \begin{split}
    \mathtt{term(dis-dec)}&=\sum_{t=1}^T\left(2\|\hat{\q}_{t}-\q_t^{(i)}\|_1+G\|\w_t^{(i)}-\w^{(j(i))}_t\|_2+G\sum_{j=1}^m q_{t,j}^{(i)}\|\w_t^{(i)}-\w^{(j)}_t\|_2\right)\\
    &\leq \sum_{t=1}^T\left(2\|\hat{\q}_{t}-\tilde{\q}_t\|_1+2\|\tilde{\q}_t-\q_t^{(i)}\|_1+4G\epswdist\right)\\
    &\leq \sum_{t=1}^T\left(2\|\hat{\q}_{t}-\tilde{\q}_t\|_1+2\sqrt{m}\|\tilde{\q}_t-\q_t^{(i)}\|_2+4G\epswdist\right)\\
    &\overset{(\ref{eqn:eps_q_1to2})}{\leq} 4\sum_{t=1}^T\left(m\epsqdist+G\epswdist\right)
    \end{split}
\end{equation}

Substituting the regret bounds~(\ref{eqn:proof:dec:regw_final_bound}) and~(\ref{eqn:proof:dec:q:final_bound}), together with the discrepancy estimate $\E[\mathtt{term(dis\text{-}dec)}]\le 4\sum_{t=1}^T \E\left[m\epsqdist+G\epswdist\right]$,
into the decomposition~(\ref{eqn:proof:dec:error_decom_final_dec}) yields, for all $i\in[m]$,
\begin{equation}\label{eqn:proof:dec:final_opt_bound_general}
\begin{split}
\E\left[\varepsilon_{\phi}(\wb_T^{(i)})\right]
\le&
\frac{1}{T}\Bigg[
\sqrt{\frac{T}{2}}DG
+\frac{DG}{\sqrt{2T}}\sum_{t=1}^T \E\left[\left(1+m\epsqdist\right)^2\right]
+3G\sum_{t=1}^T \E\left[\left(1+m\epsqdist\right)\epswdist\right] \\
&\qquad\quad
+\sqrt{m}\sum_{t=1}^T \E[\epsqdist]
+4\sqrt{T\ln m}+4\sum_{t=1}^T \E\left[m\epsqdist+G\epswdist\right]
\Bigg].
\end{split}
\end{equation}
For any $i\in[m]$ and $t>1$, Lemma~\ref{lem:Projection} and average preservation give $\|\w_t^{(i)}-\tilde{\w}_t\|_2 \le\eta_w\|\Fb_{t-1}^{(i)}-\bFb_{t-1}\|_2$. The same argument for the dual update yields $\|\q_t^{(i)}-\tilde{\q}_t\|_2\le\eta_q\|\Lb_{t-1}^{(i)}-\bLb_{t-1}\|_2$.
Thus the second-moment assumptions imply $(\E[\epswdist^2])^{1/2}\le\eta_w\varepsilon_{c,w},$ and $(\E[\epsqdist^2])^{1/2}\le\eta_q\varepsilon_{c,q}$.

Minkowski's and Cauchy--Schwarz inequalities then give
\[
\E[(1+m\epsqdist)^2]\le(1+m\eta_q\varepsilon_{c,q})^2,
\quad
\E[(1+m\epsqdist)\epswdist]
\le(1+m\eta_q\varepsilon_{c,q})\eta_w\varepsilon_{c,w}.
\]
Substituting these bounds into~(\ref{eqn:proof:dec:final_opt_bound_general}), together with
$\E[\epswdist]\le\eta_w\varepsilon_{c,w}$ and
$\E[\epsqdist]\le\eta_q\varepsilon_{c,q}$, yields, for all $i\in[m]$,
\begin{equation}\label{eqn:proof:dec:final_opt_bound:1}
\begin{split}
\E\left[\varepsilon_{\phi}(\wb_T^{(i)})\right]
\le&
\frac{DG}{\sqrt{2T}}
+\frac{DG}{\sqrt{2T}}\left(1+m^2\sqrt{\frac{2\ln m}{T}}\varepsilon_{c,q}\right)^2 
+\frac{3\sqrt{2}mD}{\sqrt{T}}
\left(1+m^2\sqrt{\frac{2\ln m}{T}}\varepsilon_{c,q}\right)\varepsilon_{c,w} \\
&+m^{3/2}\sqrt{\frac{2\ln m}{T}}\varepsilon_{c,q}+4\sqrt{\frac{\ln m}{T}}+4\left(m^2\sqrt{\frac{2\ln m}{T}}\varepsilon_{c,q}
+\frac{\sqrt{2}mD}{\sqrt{T}}\varepsilon_{c,w}\right)
\end{split}
\end{equation}

For the compact notation used in Theorem~\ref{thm:Dec_GDRO_modular}, define
\begin{equation*}\label{eqn:dec:consensus_terms}
C_1 =\left(1+m^2\sqrt{\frac{2\ln m}{T}}\,\varepsilon_{c,q}\right)^2, \ 
C_2 =m\left(1+m^2\sqrt{\frac{2\ln m}{T}}\,\varepsilon_{c,q}\right)\varepsilon_{c,w},\
C_3 =m^2\sqrt{2\ln m}\,\varepsilon_{c,q}.
\end{equation*}
Equation~(\ref{eqn:proof:dec:final_opt_bound:1}) then yields the bound in Theorem~\ref{thm:Dec_GDRO_modular}.

\subsection{The Proof of Corollary~\ref{cor:Dec_AG}}\label{app:proof:Dec_AG}
When adopting the accelerated gossip update~(\ref{lem:acc_gossip_update}), we set $\Fb_t^{(i)}=\Xb_t^{(i)}$ and $\bFb_t=\bXb_t=\frac{1}{m}\Xb_t\mathbf{1}$.
By Lemma~\ref{lem:iter_bound:ag_gossip}, the average is preserved, i.e., $\bFb_t=\bFb_{t-1}+\bg_{w,t}$, and hence the conservation condition in Theorem~\ref{thm:Dec_GDRO_modular} is satisfied.

First consider $T\ge2\ln m$. Set
$k_q=\lceil(m^{5/2}+1)^2\rceil$ and
$k_w=\lceil(m/\sqrt{\ln m}+1)^2\rceil$.
Since $\sum_i\|\g_{q,t}^{(i)}\|_2^2\le m$, Lemma~\ref{lem:iter_bound:ag_gossip} gives the uniform pathwise bound
\begin{equation}\label{eqn:proof:dec:epsqdist_AG_final}
\max_{t\in[T]}\|\Lb_t-\bLb_t\mathbf{1}^\top\|_F
\le\frac{\sqrt m}{\sqrt{k_q}-1}\le\frac1{m^2}.
\end{equation}
Consequently, $\epsqdist\le\eta_q/m^2$ for all $t$, and
\[
\sum_i\|\g_{w,t}^{(i)}\|_2^2
\le(2+2m\epsqdist^2)G^2
\le\left(2+\frac{4\ln m}{mT}\right)G^2\le4G^2.
\]
Applying Lemma~\ref{lem:iter_bound:ag_gossip} with the uniform input bound $Q=4G^2$ yields
\begin{equation}\label{eqn:proof:dec:epswdist_AG_final}
\max_{t\in[T]}\|\Fb_t-\bFb_t\mathbf{1}^\top\|_F
\le\frac{2G}{\sqrt{k_w}-1}\le\frac{2G\sqrt{\ln m}}m.
\end{equation}
Thus we may take $\varepsilon_{c,q}=1/m^2$ and
$\varepsilon_{c,w}=2G\sqrt{\ln m}/m$.
Substituting these relations into the general optimization bound~(\ref{eqn:proof:dec:final_opt_bound:1}) yields, for all $i\in[m]$,
\begin{equation}\label{eqn:proof:dec:final_bound_AG_stepsizes}
\begin{split}
\E\left[\varepsilon_{\phi}(\wb_T^{(i)})\right]
\le&
\frac{DG}{\sqrt{2T}}
+\frac{DG}{\sqrt{2T}}\left(1+\sqrt{\frac{2\ln m}{T}}\right)^2
+\frac{6\sqrt{2\ln m}DG}{\sqrt{T}}\left(1+\sqrt{\frac{2\ln m}{T}}\right) \\
&\quad
+\frac{1}{\sqrt{m}}\sqrt{\frac{2\ln m}{T}}
+4\sqrt{\frac{\ln m}{T}}
+4\sqrt{\frac{2\ln m}{T}}
+\frac{8\sqrt{2\ln m}DG}{\sqrt{T}}\\
\leq &(32DG+9\sqrt{2})\sqrt{\frac{\ln m}{T}},
\end{split}
\end{equation}
For the last inequality, set $a=\sqrt{2\ln m/T}\le1$. After factoring out $\sqrt{(\ln m)/T}$, the coefficient of $DG$ is
\[
A(m,a)=\frac{1+(1+a)^2}{\sqrt{2\ln m}}+\sqrt2(14+6a).
\]
For $m\ge3$, $A(m,a)\le20\sqrt2+5/\sqrt{2\ln3}<32$. For $m=2$, the integer horizon satisfies $T\ge2$, so $a\le\sqrt{\ln2}$ and
\[
A(2,a)\le\frac{1+(1+\sqrt{\ln2})^2}{\sqrt{2\ln2}}+\sqrt2(14+6\sqrt{\ln2})<32.
\]
The remaining coefficient is at most $9\sqrt2$, since $m^{-1/2}+2\sqrt2+4\le9$. For $T<2\ln m$, the claimed bound follows from $\varepsilon_\phi(\wb_T^{(i)})\le1$.

\subsection{The Proof of Corollary~\ref{cor:Dec_CG}}\label{app:proof:Dec_CG}
We now analyze Algorithm~\ref{alg:Dec_GDRO} under Choco-gossip.
By Lemma~\ref{lem:choco-gossip}, the average-preservation conditions hold automatically, namely $\sum_{i=1}^m \Fb_{t}^{(i)}=\sum_{i=1}^m\sum_{\tau=1}^{t}\g_{w,\tau}^{(i)}$ and $\sum_{i=1}^m \Lb_{t}^{(i)}=\sum_{i=1}^m\sum_{\tau=1}^{t}\g_{q,\tau}^{(i)}$.

First consider $T\ge2\ln m$, and set $k_q=4m^5+2$ and $k_w=m^2/(\ln m)+2$.
Lemma~\ref{lem:semi_dec:choco-gossip}, with $\sum_i\|\g_{q,t}^{(i)}\|_2^2\le m$, gives
\[
\left(\E[\|\Lb_t-\bLb_t\mathbf{1}^\top\|_F^2]\right)^{1/2}
\le2\sqrt{\frac{m}{k_q-2}}=\frac1{m^2}.
\]
Hence $\E[\epsqdist^2]\le\eta_q^2/m^4$, and
\[
\E\left[\sum_i\|\g_{w,t}^{(i)}\|_2^2\right]
\le(2+2m\E[\epsqdist^2])G^2
\le\left(2+\frac{4\ln m}{mT}\right)G^2\le4G^2.
\]
Applying Lemma~\ref{lem:semi_dec:choco-gossip} with $Q=4G^2$ therefore yields
\[
\left(\E[\|\Fb_t-\bFb_t\mathbf{1}^\top\|_F^2]\right)^{1/2}
\le\frac{4G}{\sqrt{k_w-2}}=\frac{4G\sqrt{\ln m}}m.
\]
We may thus take $\varepsilon_{c,q}=1/m^2$ and
$\varepsilon_{c,w}=4G\sqrt{\ln m}/m$ in~(\ref{eqn:proof:dec:final_opt_bound:1}), yielding
\begin{equation}\label{eqn:proof:dec:final_bound_choco_stepsizes}
\begin{split}
\E\left[\varepsilon_{\phi}(\wb_T^{(i)})\right]
\le&
\frac{DG}{\sqrt{2T}}
+\frac{DG}{\sqrt{2T}}\left(1+\sqrt{\frac{2\ln m}{T}}\right)^2
+\frac{12\sqrt{2\ln m}DG}{\sqrt{T}}\left(1+\sqrt{\frac{2\ln m}{T}}\right) \\
&\quad
+\frac{1}{\sqrt{m}}\sqrt{\frac{2\ln m}{T}}
+4\sqrt{\frac{\ln m}{T}}
+4\sqrt{\frac{2\ln m}{T}}
+\frac{16\sqrt{2\ln m}DG}{\sqrt{T}}\\
\leq& \frac{5DG}{\sqrt{2T}}+\frac{40\sqrt{2\ln m}DG}{\sqrt{T}}+9\sqrt{\frac{2\ln m}{T}} \\
\leq& (61DG+9\sqrt{2})\sqrt{\frac{\ln m}{T}}.
\end{split}
\end{equation}
where we use $T\ge2\ln m$, $m^{-1/2}+2\sqrt2+4\le9$, and $40\sqrt2+5/\sqrt{2\ln m}\le40\sqrt2+5/\sqrt{2\ln2}<61$. For $T<2\ln m$, the claimed bound follows from $\varepsilon_\phi(\wb_T^{(i)})\le1$.

\subsection{The Proof of Lemma~\ref{lem:iter_bound:ag_gossip}}
Denote $\g_t=[\g_t^{(1)},\ldots,\g_t^{(m)}]$ and $\bg_t=\frac{1}{m}\g_t\mathbf{1}$.
By the choice of $L$, we have $14 c^{2L}\le \frac{1}{k}$.
By Lemma~\ref{lem:ag_convergence} and~(\ref{lem:acc_gossip_update}), Algorithm~\ref{alg:AG} preserves the average, i.e.,
\[
\bXb_t
=\frac{1}{m}\Xb_t\mathbf{1}
=\frac{1}{m}(\Xb_{t-1}+\g_t)\mathbf{1}
=\bXb_{t-1}+\bg_t.
\]
Therefore,
\begin{align*}
\|\Xb_t-\bXb_t\mathbf{1}^\top\|_F
&\le \sqrt{14} c^{L}
\|\Xb_{t-1}+\g_t-(\bXb_{t-1}+\bg_t)\mathbf{1}^\top\|_F \\
&\le \frac{1}{\sqrt{k}}
\|(\Xb_{t-1}-\bXb_{t-1}\mathbf{1}^\top)+(\g_t-\bg_t\mathbf{1}^\top)\|_F \\
&\le \frac{1}{\sqrt{k}}\|\Xb_{t-1}-\bXb_{t-1}\mathbf{1}^\top\|_F
   +\frac{1}{\sqrt{k}}\|\g_t-\bg_t\mathbf{1}^\top\|_F,
\end{align*}

Since $\bg_t$ is the minimizer of $\sum_i\|\g_t^{(i)}-\u\|_2^2$ over $\u$, we have the standard bound
\[
\|\g_t-\bg_t\mathbf{1}^\top\|_F^2
=\sum_{i=1}^m\|\g_t^{(i)}-\bg_t\|_2^2
\le \sum_{i=1}^m\|\g_t^{(i)}\|_2^2
\le Q.
\]
Define $\e_t=\|\Xb_t-\bXb_t\mathbf{1}^\top\|_F$. Then we obtain $\e_t \le \frac{1}{\sqrt{k}}\e_{t-1}+\frac{\sqrt{Q}}{\sqrt{k}}$.
By recursion and $\e_0=0$, we have
\begin{align*}
\e_t
&\le \frac{\sqrt{Q}}{\sqrt{k}}\sum_{s=0}^{t-1}\left(\frac{1}{\sqrt{k}}\right)^s \leq \frac{\sqrt{Q}}{\sqrt{k}-1},
\end{align*}
which proves~(\ref{eq:lem:ag_gossip}).

\subsection{The Proof of Lemma~\ref{lem:semi_dec:choco-gossip}}
Recall that the number of communication rounds is defined as $L = \left\lceil \frac{2 \ln k}{\gamma \rho} \right\rceil$ for a constant \( k > 2 \). Then denote $\g_t=[\g_{t}^{(1)},\ldots,\g_{t}^{(m)}]$ and $\bg_{t}=\frac{1}{m}\g_t\mathbf1$ we have $\left(1 - \frac{\gamma \rho}{2} \right)^{L} \leq \frac{1}{k}$.
Let
$a_t=\E[\|\Xb_t-\bXb_t\mathbf{1}^\top\|_F^2+\|\Xb_t-\hXb_t\|_F^2]$.
Conditioning on the inputs to the current consensus call, applying
Lemma~\ref{lem:choco-gossip}, and then taking total expectations gives
\begin{align*}
a_t
&\le\frac1k\E\left[
\|\Xb_{t-1}+\g_t-(\bXb_{t-1}+\bg_t)\mathbf{1}^\top\|_F^2
+\|\Xb_{t-1}+\g_t-\hXb_{t-1}\|_F^2\right]\\
&\le\frac2k a_{t-1}
+\frac2k\E\left[\sum_{i=1}^m
\left(\|\g_t^{(i)}-\bg_t\|_2^2+\|\g_t^{(i)}\|_2^2\right)\right]\\
&\le\frac2k a_{t-1}+\frac{4Q}{k},
\end{align*}
where the last inequality uses
$\sum_i\|\g_t^{(i)}-\bg_t\|_2^2\le\sum_i\|\g_t^{(i)}\|_2^2$.
Since $a_0=0$, unrolling the recursion gives
\[
a_t\le\frac{4Q}{k}\sum_{s=0}^{t-1}(2/k)^s
\le\frac{4Q}{k-2},
\]
which proves~(\ref{eqn:lem:choco-gossip}).

\end{document}